\documentclass[reqno]{amsart}
\theoremstyle{plain}
\usepackage{latexsym}
\usepackage{amssymb, color}
\usepackage{enumerate}
\usepackage{array}
\usepackage{lipsum}
\usepackage{epsfig}
\usepackage[sc]{mathpazo}
\def\Xint#1{\mathchoice
	{\XXint\displaystyle\textstyle{#1}}%
	{\XXint\textstyle\scriptstyle{#1}}%
	{\XXint\scriptstyle\scriptscriptstyle{#1}}%
	{\XXint\scriptscriptstyle\scriptscriptstyle{#1}}%
	\!\int}

\def\XXint#1#2#3{{\setbox0=\hbox{$#1{#2#3}{\int}$}
		\vcenter{\hbox{$#2#3$}}\kern-.5\wd0}}

\newcommand{\Cb}{\mathbb{C}}

\newcommand{\Rb}{\mathbb{R}}

\newcommand{\lt}{\left}
\newcommand{\rt}{\right}
\newcommand{\nl}{\newline}

\newcommand{\nn}{\nonumber}

\newcommand{\wii}{\varpi}

\newcommand{\qd}{\quad}
\newcommand{\spt}{\mathrm{Spt}}

\newcommand{\ep}{\epsilon}

\newcommand{\wt}{\widetilde}

\newcommand{\OI}{\mathcal{O}}

\newcommand{\na}{\nabla}

\newcommand{\R}{\mathrm {I\!R}}

\newcommand{\dia}{\diamondsuit}

\newcommand{\red}[1]{{\textcolor{black}{#1}}} 
\newcommand{\blue}[1]{{\textcolor{black}{#1}}} 
\newcommand{\xred}[1]{{\textcolor{black}{#1}}} 
\newcommand{\xblue}[1]{{\textcolor{black}{#1}}} 
\newcommand{\bblue}[1]{{\textcolor{black}{#1}}} 
\newcommand{\green}[1]{{\textcolor{black}{#1}}}

\newtheorem{a1}{Lemma}
\newtheorem{a2}[a1]{Theorem}

\newtheorem{a6}[a1]{Corollary}

\theoremstyle{remark}

\begin{document}
	
	\title[Regularity of the Eikonal equation with two vanishing entropies ]
	{Regularity of the Eikonal equation with two vanishing entropies  }
	\author{Andrew Lorent, Guanying Peng}
	\address{Mathematics Department\\University of Cincinnati\\2600 Clifton Ave.\\ Cincinnati OH 45221 }
	\email{lorentaw@uc.edu, penggg@ucmail.uc.edu}
	\subjclass[2000]{28A75}
	\keywords{eikonal equation, Aviles Giga functional,  entropies, non-linear Beltrami equation, differential inclusions}
	\maketitle

	\begin{abstract}
		
		The Aviles-Giga functional $I_{\ep}(u)=\int_{\Omega} \frac{\lt|1-\lt|\na u\rt|^2\rt|^2}{\ep}+\ep \lt|\na^2 u\rt|^2 \; dx$ is a well known second order functional that models phenomena from blistering to liquid crystals. The zero energy states of the Aviles-Giga functional have been characterized by 
		Jabin, Otto, Perthame \cite{otto}. Among other results they showed that if $\lim_{n\rightarrow \infty} I_{\ep_n}(u_n)=0$ for some sequence $u_n\in W^{2,2}_0(\Omega)$ and $u=\lim_{n\rightarrow \infty} u_n$ then $\na u$ is Lipschitz continuous outside a finite set. This is essentially a corollary to their theorem that if $u$ is a solution to the Eikonal equation $\lt|\na u\rt|=1$ a.e.\ and if for every "entropy"  $\Phi$ (in the sense of \cite{mul2}, Definition 1) function $u$ satisfies $\nabla\cdot\lt[\Phi(\na u^{\perp})\rt]=0$ distributionally in $\Omega$ then $\na u$ is \bblue{locally} Lipschitz continuous outside a \bblue{locally} finite set.
		
		In this paper we generalize this result by showing \bblue{that} if $u$ satisfies the Eikonal equation and if 
		\begin{equation}
		\label{eqi88}
		\nabla\cdot\lt(\bblue{\wt\Sigma_{e_1 e_2}}(\na u^{\perp})\rt)=0\text{ and }\nabla\cdot\lt(\bblue{\wt\Sigma_{\ep_1 \ep_2}}(\na u^{\perp})\rt)=0\text{ distributionally in }\Omega, 
		\end{equation}
		where  $\wt{\Sigma}_{e_1 e_2}$ and $\wt{\Sigma}_{\ep_1 \ep_2}$ are the entropies introduced by Ambrosio, DeLellis, Mantegazza \cite{ADM}, Jin, Kohn \cite{kohn}, then  $\na u$ is \bblue{locally} Lipschitz continuous outside a \bblue{locally} finite set. Condition (\ref{eqi88}) being fairly natural this result could also be considered a contribution to the study of the regularity of solutions of the Eikonal equation. The method of proof is to transform any solution of the Eikonal equation satisfying (\ref{eqi88}) into a differential inclusion $DF\in K$ where $K\subset M^{2\times 2}$ is a connected compact set of matrices without Rank-$1$ connections. Equivalently this differential inclusion can be written as \xred{a constrained} non-linear Beltrami equation. The set $K$ is also non-elliptic in the sense of Sverak \cite{sverak}. By use \bblue{of} this transformation and by utilizing ideas from the work on regularity of solutions of the Eikonal equation in fractional Sobolev space by Ignat \cite{ignat},  DeLellis, Ignat \cite{DeI} as well as methods of Sverak \cite{sverak}, regularity is established.  
	\end{abstract}

	\section{Introduction}
	\label{intro}
	
	The Aviles-Giga functional is the second order functional
	\begin{equation*}
	I_{\ep}(u)=\int_{\Omega} \frac{\lt|1-\lt|\na u\rt|^2\rt|^2}{\ep}+\ep \lt|\na^2 u\rt|^2 \; dx
	\end{equation*}
	minimized over the space of functions $W_0^{2,2}(\Omega;\R)$ or $W_0^{2,2}(\Omega;\R)\cap \lt\{u: \na u(x)=\eta_x\text{ on }\partial \Omega\rt\}$ where $\eta_x$ is the inward pointing unit normal to $\partial \Omega$. The Aviles-Giga functional $I_{\ep}$ forms a model for blistering  and (in
	certain regimes) a model for liquid crystals \cite{avgig0}, \cite{kohn}, \cite{gioort}. In addition there is a closely related functional modeling thin magnetic
	films known as the \em micromagnetics functional \rm \cite{mul2}, \cite{deott1}, \cite{cdeott}, \cite{ser1}, \cite{ser3}, \cite{ser2}, \cite{amleri}. For function $u\in W_0^{2,2}(\Omega)$ we 
	refer to $I_{\ep}(u)$ as the \em Aviles-Giga energy \rm of $u$. The Aviles-Giga functional is the most natural higher order generalization of the Modica-\bblue{Mortola} functional \cite{modmor}.
	
	The biggest open problem in the study of the Aviles-Giga functional is the characterization of \xred{its $\Gamma$-limit, \cite{avgig0}, \cite{avgig1}, \cite{kohn}, \cite{ADM}}.  Given the structure of $I_{\ep}$ it is not a surprise that the conjectured limiting function class is a subspace of functions that satisfy the Eikonal equation 
	\begin{equation}
	\label{eqi4}
	\lt|\na u(x)\rt|=1\text{ for }a.e.\ x\in \Omega.
	\end{equation}  
	By analogy to the Modica-Mortola functional, it might be expected that the limiting function space is also a subspace of $\lt\{v: \na v\in BV \rt\}$ and the limiting energy is related to $\|D \lt[\na u\rt]\|$. However this is completely \em false\rm; see the example following Theorem 3.9 of \cite{ADM}. It is necessary to build a function class that is in a sense analogous to the function class  $\lt\{v: \na v\in BV \rt\}$ that is \bblue{tailored} to the functional $I_{\ep}$. This is done by introducing certain \em entropies \rm on the space of solutions of the Eikonal equation. The divergence of these entropies will (by virtue of the structure of $I_{\ep}$) form measures that in regular examples pick up the jump in the gradient $\na u$. Specially it can be shown \cite{ADM}, \cite{mul2} that if $u_n\in W^{2,2}_0(\Omega)$  with the 
	property that $\limsup_{n\rightarrow \infty} I_{\ep_n}(u_n)<\infty$ then for some subsequence \bblue{$\{n_k\}$} we have $u_{\bblue{n_k}}\overset{W^{1,3}(\Omega)}{\rightarrow} u$. This allows us to show that if the vector field  $\Sigma_{\xi \eta}u$ is defined by 
	\begin{equation*}
	%\label{eqzzz221}
	\Sigma_{\xi\eta}u:=u_{\xi}\lt(1-u_{\eta}^2-\frac{1}{3}u_{\xi}^2\rt)\xi-u_{\eta}\lt(1-u_{\xi}^2-\frac{1}{3}u_{\eta}^2\rt)\eta,
	\end{equation*}
	(where $u_{\xi}$ and $u_{\eta}$ are the partial derivatives along $\xi$ and $\eta$ respectively) then $\nabla\cdot\lt(\Sigma_{\xi \eta}u \rt)$ is a measure. 
	So instead of having that the gradient of the gradient is a measure (as would be the case if 
	$u\in \lt\{v: \na v\in BV \rt\}$) we have that the divergence of a vector field made up of first order partial gradients is a measure, which ''morally''  is not that far away.  
	
	\bblue{Following \cite{ADM},} we denote by $(e_1,e_2)$ the canonical basis of $\Rb^2$, and by 
	\begin{equation}
	\label{eqzzz222}
	\ep_1:=\lt(\frac{1}{\sqrt{2}},\frac{1}{\sqrt{2}}\rt), \quad\ep_2=\lt(-\frac{1}{\sqrt{2}},\frac{1}{\sqrt{2}}\rt)
	\end{equation}
	the basis obtained from $(e_1,e_2)$ under an anticlockwise rotation of $\frac{\pi}{4}$. It is \bblue{straightforward to check} that 
	\begin{equation}
	\label{eqi501}
	\Sigma_{e_1 e_2}u=\lt(u_{,1}\lt(1-u_{,2}^2-\frac{u_{,1}^2}{3}\rt),  -u_{,2}\lt(1-u_{,1}^2-\frac{u_{,2}^2}{3}\rt)\rt)
	\end{equation}
	and
	\begin{equation}
	\label{eqi502}
	\Sigma_{\ep_1 \ep_2}u=\lt(u_{,2}\lt(1-\frac{2 u_{,2}^2}{3}\rt),  u_{,1}\lt(1-\frac{2 u_{,1}^2}{3}\rt)\rt).
	\end{equation}
	It has been shown in \cite{ADM} that the measure 
	\begin{equation*}
	%\label{eqi6}
	S\rightarrow \lt|\lt|\lt(\begin{array}{c}  \nabla\cdot\lt(\Sigma_{e_1 e_2}u  \rt) \\
	\nabla\cdot\lt(\Sigma_{\ep_1 \ep_2}u\rt) \end{array}\rt)\rt|\rt|(S)\text{ for any }S\subset \R^2
	\end{equation*}
	forms a lower bound on the energy $I_{\ep_n}(u_{n})$ of any sequence \bblue{$\{u_n\}$} such that $\lim_{n \rightarrow \infty} u_n=u$. As such the functional 
	\begin{equation}
	\label{eqi7}
	u\rightarrow \lt|\lt|\lt(\begin{array}{c}  \nabla\cdot\lt(\Sigma_{e_1 e_2}u  \rt) \\
	\nabla\cdot\lt(\Sigma_{\ep_1 \ep_2}u\rt) \end{array}\rt)\rt|\rt|(\Omega)
	\end{equation}
	was conjectured in \cite{ADM} to be the $\Gamma$-limiting energy of the Aviles-Giga functional. 
	
	Following \cite{mul2}, \cite{ottodel1} we say $\Phi\in C_{c}^{\bblue{\infty}}(\R^2;\R^2)$ is an entropy if 
	\begin{equation}
	\label{eqi151}
	z\cdot D\Phi(z) z^{\perp}=0\text{ for all }z\bblue{\in\R^2}, \Phi(0)=0, D\Phi(0)=0, 
	\end{equation}
	\bblue{where $z^{\perp}=(-z_2,z_1)$ is the anticlockwise rotation of $z$ by $\frac{\pi}{2}$.} Vector fields 
	\begin{equation}
	\label{eqi551}
	\wt{\Sigma}_{e_1 e_2}(x,y):=\lt( y\lt(1-x^2-\frac{y^2}{3}\rt),x\lt(1-y^2-\frac{x^2}{3}\rt)\rt)\text{ and  }\wt{\Sigma}_{\ep_1 \ep_2}(x,y):=\lt( -x\lt(1-\frac{2 x^2}{3}\rt),y\lt(1-\frac{2 y^2}{3}\rt)\rt)
	\end{equation}
	satisfy 
	\begin{equation}
	\label{eqi301}
	z\cdot D \wt{\Sigma}_{e_1 e_2}(z) z^{\perp}=0\text{ for all }z\in \mathbb{S}^1, \;\; \wt{\Sigma}_{e_1 e_2}(0)=0
	\end{equation} 
	and 
	\begin{equation}
	\label{eqi302}
	z\cdot D \wt{\Sigma}_{\ep_1 \ep_2}(z) z^{\perp}=0\text{ for all }z\in \mathbb{S}^1,\;\;  \wt{\Sigma}_{\ep_1 \ep_2}(0)=0.
	\end{equation}  
	Note that $\Sigma_{e_1 e_2} u\overset{(\ref{eqi501}),(\ref{eqi551})}{=}\wt{\Sigma}_{e_1 e_2}(\na u^{\perp})$ and  $\Sigma_{\ep_1 \ep_2} u\overset{(\ref{eqi502}),(\ref{eqi551})}{=}\wt{\Sigma}_{\ep_1 \ep_2}(\na u^{\perp})$, \bblue{where $\na u^{\perp}=(-u_{,2},u_{,1})$}. Since we are applying $\wt{\Sigma}_{e_1 e_2}$, $\wt{\Sigma}_{\ep_1 \ep_2}$ to gradient vector fields $\na u$ that satisfy $\lt|\na u\rt|=1$ a.e., for simplicity, and following the convention of \cite{ottodel1}, we will call them entropies even though they only satisfy (\ref{eqi301}), (\ref{eqi302}). However this is just a naming convenience \bblue{and} is not important to the mathematics that follows. Whenever we use any results about entropies from \cite{mul2} we will mean vector fields  $\Phi\in C_c^{\bblue{\infty}}(\R^2;\R^2)$ that satisfy (\ref{eqi151}). The main point about entropies is that given a sequence $\{u_n\}$ that satisfies  
	$\limsup_{n\rightarrow \infty} I_{\ep_n}(u_n)<\infty$ and $u=\lim_{n\rightarrow \infty} u_n$, if $\Phi$ is an entropy then 
	$\nabla\cdot \lt[\Phi\lt(\na u^{\perp}\rt)\rt]$ is a measure.  
	
	The characterization of this class of entropies is one of the main \bblue{achievements} of \cite{mul2} and it \bblue{leads} to many further developments. It was the main tool used in \cite{mul2} to prove pre-compactness in $W^{1,3}(\Omega)$ of a sequence of functions $\{u_n\}$ of bounded Aviles-Giga energy (an alternative proof just using \em two \rm entropies $\wt{\Sigma}_{e_1 e_2}$,  $\wt{\Sigma}_{\ep_1 \ep_2}$ is \bblue{provided} in \cite{ADM}). More importantly it allows for the classification achieved by Jabin, Otto, Perthame  in \cite{otto} of all functions $u$ and all domains $\Omega$ for which there exists a sequence $\bblue{\{u_n\}\subset} W^{2,2}_0(\Omega)$ such that $u=\lim_{n\rightarrow \infty} u_n$ and $\lim_{n\rightarrow \infty} I_{\ep_n}(u_n)=0$. Functions $u$ with this property are called zero energy states. It was shown in \cite{otto} that if $\Omega\not=\R^2$ then $\Omega$ is a ball and (after possibly change of sign) $u$ is just the distance function away from the \green{boundary of the} ball. The characterization of entropies also permitted the deep work on the structure of solutions of the Eikonal equation $u$ that arise as limits of sequences of \em finite \rm Aviles-Giga energy \cite{ottodel2}. 
	
	While the works \cite{ottodel2}, \cite{otto} are impressive achievements and indeed represent the state of the art with respect to the structure of solutions of the Eikonal equation that arise as limits of sequences of finite (or converging to zero) Aviles-Giga energy, when these \bblue{results} are formulated simply in terms of the Eikonal equation, the statements can appear a bit technical. 
	\begin{a2}[\bblue{\cite{otto}}]
		\label{JOP1}
		Let $\Omega$ be any open set in $\R^2$. Let $m\bblue{:\Omega\rightarrow \R^2}$ be a measurable function that satisfies $\lt|m(x)\rt|=1$ for a.e. $x\in \Omega$ and 
		\begin{equation}
		\label{eqi45}
		\xi \cdot \na \chi\lt(\cdot,\xi \rt)=0\text{ distributionally in }\Omega\text{ for all }\xi\in \mathbb{S}^1,
		\end{equation}
		where 
		\begin{equation}
		\label{eqi45.6}
		\chi(x,\xi):=\begin{cases}
		1 & \text{ for }m(x)\cdot \xi>0,\\
		0 & \text{ for }m(x)\cdot \xi\leq 0. 
		\end{cases}
		\end{equation}	
		Then $m$ is locally Lipschitz outside a \bblue{locally} finite set of points. 
	\end{a2} 
	It turns out that $\xi\chi(\cdot,\xi)$ is the pointwise limit of a sequence of entropies $\lt\{\Phi_n\rt\}$ (see the proof of \bblue{Lemma 4, \cite{mul2}}), so if vector field $m$ is such that 
	\begin{equation}
	\label{eqi998}
	\nabla\cdot\lt[\Phi\lt(m\rt)\rt]=0\text{ distributionally in }\Omega \text{ for all entropies }\Phi,
	\end{equation}
	then $m$ satisfies (\ref{eqi45}). Hence by Theorem \ref{JOP1} any vector field $m$ satisfying (\ref{eqi998})  is locally Lipschitz outside a \bblue{locally} finite set of points. This is the main result needed by Jabin, Otto, Perthame \cite{otto} to characterize all zero energy states of the Aviles-Giga energy. 
	
	\begin{a6}[\bblue{\cite{otto}}] 
		\label{JOP2}
		Let $u$ be a limit of a sequence $\lt\{u_n\rt\}\subset W^{2,2}_0(\Omega)$ with $\lim_{n\rightarrow \infty} I_{\ep_n}(u_n)=0$ then $\na u$ is Lipschitz outside a finite set of points.
	\end{a6}
	Actually in \cite{otto} a more general result is proved that includes zero energy states of the micromagnetic \bblue{functional}, but since our interest \bblue{is} focused on the Aviles-Giga \bblue{functional} we do not state their result in full generality.  
	
	What is achieved in this paper is a proof of the regularity result under the much weaker condition that only the divergence of $\na u^{\perp}$ applied to two entropies $\wt{\Sigma}_{e_1 e_2}$ and $\wt{\Sigma}_{\ep_1 \ep_2}$ vanishes.  
	\begin{a2} 
		\label{T1}
		Let $\Omega \subset \R^2$ be \bblue{a bounded simply-connected domain} and $u$ be a solution to the Eikonal equation (\ref{eqi4}). Suppose
		\begin{equation}
		\label{eqi46}
		\nabla\cdot\lt(\Sigma_{e_1 e_2}u  \rt)=0\text{ and }\nabla\cdot\lt(\Sigma_{\ep_1 \ep_2}u  \rt)=0\text{ distributionally in }\Omega. 
		\end{equation}
		Then $\na u$ is locally Lipschitz outside a \bblue{locally} finite set of points $S$. Moreover, in any convex neighborhood $\OI\subset \subset \Omega$ of a point $\zeta\in S$ there exists $\bblue{\alpha}\in \lt\{-1,1\rt\}$ such that 
		\begin{equation}
		\label{eq210}
		u(x)=\bblue{\alpha} \lt|x-\zeta\rt|\text{ for any }x\in \OI. 
		\end{equation}
	\end{a2}
	
	This result also includes Corollary \ref{JOP2} as a consequence \bblue{in the case that $\Omega$ satisfies the assumptions of Theorem \ref{T1}}. The value of this result is twofold. Firstly the Eikonal equation is a much studied equation whose more general form $\lt|\na u\rt|=f$ occurs in numerous areas of physics (geometric optics, wave propagation) and applied mathematics. Historically there has been great interest in first uniqueness and then subsequently regularity of the 
	Eikonal equation. Uniqueness was largely resolved by the development of the regularity of viscosity solutions \cite{crandal}, \green{and} subsequent regularity results have been established by a number of authors, \green{\cite{can}, \cite{can2}}. \bblue{Indeed, regularity} and uniqueness of the Eikonal equation was one of the early triumphs that follwed the development of the theory of viscosity solutions. The Eikonal equation with the additional assumption of two vanishing entropies seems to us a fairly natural condition and as such the statement of Theorem \ref{T1} is of interest purely from the perspective of the Eikonal equation alone. On this topic we mention the recent powerful results of Ignat \cite{ignat} and DeLellis, Ignat \cite{DeI} on regularity of solutions of the Eikonal equation in fractional Sobolev spaces. We learned a great deal and took numerous ideas from these works. 
	
	Our principle interest however is in the Aviles-Giga functional. As previously described the original conjectured $\Gamma$-limiting energy from \cite{ADM} is given by (\ref{eqi7}). As the study of the Aviles-Giga functional evolved it was increasingly understood that to make progress the conjectured $\Gamma$-limiting energy had to be an energy that incorporated all the entropies, not simply $\wt{\Sigma}_{e_1 e_1}$ and $\wt{\Sigma}_{\ep_1 \ep_2}$. As mentioned the proof of Corollary \ref{JOP2} requires the use of a sequence of entropies $\lt\{ \Phi_n \rt\}$ that approximates $\bblue{\xi\chi(\cdot,\xi)}$. In \cite{ottodel1} DeLellis, Otto proved many strong structural results on a class of solutions of the Eikonial equation \bblue{denoted by} $A(\Omega)$ that includes all $W^{1,3}(\Omega)$ limits of sequences $\bblue{\{u_n\}\subset} W^{2,2}_0(\Omega)$ that have equibounded Aviles-Giga energy. Among the results they proved was that for any $u\in A(\Omega)$ there exists a set of $\sigma$-finite $H^1$ measure $J$ on which $\bblue{\na}u$ has jumps and has traces in exactly the way it would have if $\bblue{\na}u\in BV$. What would be most natural is if $J$ was the singular set of \xred{vector valued} measure that is the $\Gamma$-limiting energy of $I_{\ep}$. However this is not exactly the case and $J$ has to be defined as the singular set of measure into the dual space of all entropies (see \cite{ottodel1}, proof of Proposition \bblue{1}). It is in some sense a singular set of an infinite set of entropies simultaneously.
	
	While utilizing the information available from all entropies is in our opinion the best way to progress with the study of the Aviles-Giga functional, it does have the disadvantage that the statements of the theorems proved are less transparent. It is for example not clear what the conjecture for the $\Gamma$-limiting energy of the Aviles-Giga energy is. What Theorem \ref{T1} does is \bblue{to} raise the possibility of reformulating the structure results of \cite{ottodel1}, \cite{otto} in terms of the two entropies $\wt{\Sigma}_{e_1 e_2}$, $\wt{\Sigma}_{\ep_1 \ep_2}$. Were this to be accomplished it would return the measure $S\rightarrow   \lt|\lt|\lt(\begin{array}{c}  \nabla\cdot\lt(\Sigma_{e_1 e_2}u  \rt) \\
	\nabla\cdot\lt(\Sigma_{\ep_1 \ep_2}u\rt) \end{array}\rt)\rt|\rt|(S)$ as the natural conjecture for $\Gamma$-limiting energy for the Aviles-Giga functional. 
	
	\subsection{Reduction to differential inclusions} \bblue{We denote
	\begin{equation}
	\label{eqz4}
	E(\Omega):=\{u\in W^{1,\infty}(\Omega): |\na u|=1 \text{ a.e. and (\ref{eqi46}) is satisfied}\}.
	\end{equation}}
	The starting point for our work is the transformation of functions $u\in E(\Omega)$ into functions $F_u:\Omega\rightarrow \R^2$ that satisfy the  differential inclusions $DF_u\in K$, where $K\subset M^{2\times 2}$ is a compact connected set defined by (\ref{eqy29}). This can be done because (\ref{eqi46}) can be rewritten as 
	\begin{equation*}
	%\label{eqi46.4}
	\mathrm{curl}\lt((\Sigma_{e_1 e_2}u)^{\perp}  \rt)=0\text{ and }\mathrm{curl}\lt((\Sigma_{\ep_1 \ep_2}u)^{\perp}  \rt)=0\text{ distributionally in }\Omega. 
	\end{equation*}
	Hence we can find some potential $F^1_u$ such that $\na F^1_u=(\Sigma_{e_1 e_2}u)^{\perp}$ and $F^2_u$ such that $\na F^2_u=(\Sigma_{\ep_1 \ep_2}u)^{\perp}$.\footnote{The idea to study $F_u$ comes from \cite{ADM}, see the proof of Proposition 4.6.} The structure of  $\Sigma_{e_1 e_2}$, $\Sigma_{\ep_1 \ep_2}$ implies that $DF_u\in K$ a.e. \bblue{in $\Omega$.} It is a calculation to see that $K$ does not have rank-$1$ connections, i.e., there \bblue{do} not exist $A,B\in K$, $A\not=B$ with $\mathrm{Rank}(A-B)=1$. Regularity of differential inclusions into sets without rank-$1$ connections has been studied by Sverak in his seminar paper \cite{sverak}. He showed that if function $v$ satisfies $Dv\in S\subset M^{2\times 2}$ where $S$ has no rank-$1$ connections and is elliptic in the sense that if $A,B\in S$, then $\det(A-B)\geq c\lt|A-B\rt|^2$, then $v$ is smooth. The set $K$ defined by (\ref{eqy29}) is   \em not elliptic \rm in the sense of Sverak, but it turns out 
	that for some constant $c_0\bblue{>0}$,  $\det(A-B)\geq c_0\lt|A-B\rt|^4$ for any $A,B\in K$. This is not enough to establish smoothness of $F_u$ (indeed since 
	$\na u^{\perp}$ could be a vortex, smoothness of $F_u$ could not be true) by using the methods of \cite{sverak}, \bblue{but is} enough to establish fractional Sobolev regularity. The differential inclusion $D\bblue{F_u}\in K$ can be reformulated as a \bblue{constrained} non-linear Beltrami equation and our proof of fractional Sobolev regularity can hence be formulated as the following theorem. 
	
	\begin{a2}
		\label{T2}
		Given \bblue{a bounded simply-connected domain} $\wt{\Omega}\subset \mathbb{C}$, and $v\in W^{1,\infty}(\wt{\Omega};\mathbb{C})$ that satisfies \bblue{$v(0)=0$ (assuming $0\in\wt\Omega$) and} the non-linear Beltrami system 
		\begin{equation}
		\label{eqy60}
		\frac{\partial v}{\partial \bar{z}}(z)=\blue{\frac{4}{3}}\lt(\frac{\partial v}{\partial z}(z) \rt)^3, \lt|\frac{\partial v}{\partial z}(z)\rt|=\frac{1}{2}\qd\text{ for \bblue{a.e.} }z\in \wt{\Omega},
		\end{equation}
		\bblue{we have that} 
		\begin{equation}
		\label{eqy60.4}
		\red{Dv\in W^{\bblue{\sigma},4}_{\bblue{loc}}(\wt\Omega)\qd\text{ for all } \bblue{\sigma}\in \lt(0,\frac{1}{3}\rt)}. 
		\end{equation}
		
		In addition, \bblue{given} $\wt{\Omega}'\subset \subset \wt{\Omega}$, \bblue{for all} $\ep\in \lt(0,\bblue{\frac{1}{2}}\mathrm{dist}(\wt{\Omega}',\partial \wt{\Omega})\rt)$, we have that 
		\begin{equation}
		\label{eqss1}
		\int_{\wt{\Omega}'}\int_{B_{\ep}\bblue{(0)}} \frac{\lt|Dv(z+y)-Dv(z)\rt|^4}{\ep^{2+\frac{4}{3}}} dy dz<C
		\end{equation}
		for some constant $C$ independent of $\ep$.
	\end{a2}
	
	Now if we define $\mathcal{H}_0\lt(\xi\rt):=\bblue{\frac{4}{3}}\xi^3$ then (\ref{eqy60}) can be written as $\frac{\partial v}{\partial \bar{z}}(z)=\mathcal{H}_0\lt(\frac{\partial v}{\partial z}(z)\rt)$, $\lt|\frac{\partial v}{\partial z} \rt|=\frac{1}{2}$. We will call this a \xred{contrained} non-linear Beltrami \xred{equation}. The study of equations of the form  $\frac{\partial v}{\partial \bar{z}}=\mathcal{H}\lt(z,\bblue{\frac{\partial v}{\partial z}}\rt)$  has flourished in the last few years.  Under the assumptions that 
	\begin{enumerate}[(I)]
		\item  $z\rightarrow \mathcal{H}\lt(z,w\rt)$ is measurable 
		\item  And for $w_1,w_2\in \mathbb{C}$, $\lt|\mathcal{H}\lt(z,w_1\rt)-\mathcal{H}\lt(z,w_2\rt)\rt|\leq k\lt|w_1-w_2\rt|$ for some $k<1$
	\end{enumerate}
	the existence and regularity theory of non-linear Beltrami equations resembles that of the linear theory; see \cite{boj1}, \cite{iwan}, \cite{boj2}, \cite{ast}, \cite{ast2}. But note  when restricted to the circle $\partial B_{\frac{1}{2}}(0)$ the Lipschitz constant \xred{of $\mathcal{H}_0$ }is exactly $1$, so in some sense $\mathcal{H}_0$ is a critical case \footnote{If instead we had $\mathcal{H}_0\lt(\xi\rt)=\frac{4}{3}\xi^3$ and $\lt|\frac{\partial v}{\partial z} \rt|=\alpha$ for some $\alpha\in (0,\frac{1}{2})$ we believe the standard methods of \cite{ast2} would give regularity.}. We are not aware of any other regularity results for non-linear Beltrami equations without the assumptions (I), (II). While Theorem \ref{T2} is essentially a regularity result for differential inclusions, we formulate it in the language of non-linear Beltrami equations because these are much better \bblue{known} and more studied objects. We also find the connection to this area is interesting and potentially worth further investigation.

	The connection between Theorem \ref{T1} and Theorem \ref{T2} is made by the following result.
	
	\begin{a2} 
		\label{T3}
		Let $\Omega \subset \R^2$ be \bblue{a bounded simply-connected domain}. Define $\wt{\Omega}:=\lt\{x_1+ix_2\in\Cb:(x_1,x_2)\in \Omega\rt\}$ and define $B(\wt{\Omega})$ as the set of functions $v \in W^{1,\infty}(\wt{\Omega};\mathbb{C})$ that satisfy \bblue{$v(0)=0$ (assuming $0\in\Omega$) and} the \xred{contrained} non-linear Beltrami equation 
		\begin{equation}
		\label{eqs6.5}
		\frac{\partial v}{\partial \bar{z}}(z)=\blue{\frac{4}{3}}\lt(\frac{\partial v}{\partial z}(z) \rt)^3, \lt|\frac{\partial v}{\partial z}(z)\rt|=\frac{1}{2}\qd\text{ for \bblue{a.e.} }z\in \wt{\Omega}.
		\end{equation}
		Then there exists an injective transformation 
		$$
		\Gamma:\blue{\lt[E(\Omega)/\R\rt]}\rightarrow B(\wt{\Omega}),
		$$
		where $E(\Omega)$ is defined in (\ref{eqz4}) \blue{and two functions $u_1,u_2 \in E(\Omega)$ satisfy \bblue{$u_1=u_2$ in $\lt[E(\Omega)/\R\rt]$} if and only if $u_1=u_2+C$ for some constant $C$.} Further $\Gamma$ restricted to \red{$\lt[E(\Omega)/\R\rt]\cap W^{2,1}(\Omega)$ forms a bijective transformation onto $B(\wt{\Omega})\cap W^{2,1}(\wt{\Omega})$}.
	\end{a2}
	
	However Theorem \ref{T2} and Theorem \ref{T3} will only give us fractional Sobolev \bblue{regularity}.  Ignat \cite{ignat} studied regularity of solutions of the Eikonal equation in fractional Sobolev space, \bblue{and} showed that if $u$ is a solution of the Eikonal equation and $\na u\in W^{\frac{1}{p},p}_{\bblue{loc}}\lt(\Omega\rt)$ for some $p\in \lt[1,2\rt]$ then $\na u$ is \bblue{locally} Lipschitz outside a \bblue{locally} finite set of points. Note \bblue{that} if $\na u$ \bblue{is} smooth and $\Phi$ is an entropy it follows from properties of entropies from \cite{mul2} (see Lemma \ref{l18} \bblue{of this paper}) that $\nabla\cdot\lt[\Phi\lt(\na u^{\perp}\rt)\rt]=0$. The proof of \cite{ignat} carefully exploits the structure of entropies to weaken the hypothesis on $\na u$ to that of fractional Sobolev space.  Following this work DeLellis and Ignat \cite{DeI} substantially weakened the hypothesis to $\na u\in W^{\frac{1}{p},p}_{\bblue{loc}}\lt(\Omega\rt)$ for some $p\in \lt[1,3\rt]$. It again was achieved by very careful work using the structure of entropies and by use of an estimate of Constantin, E and Titi \cite{titi}. However close though it is, this result is not quite what we need because it requires a full $1/3$ of a derivative and with the methods of \cite{sverak}, a $1/3$ of a derivative is not available - Theorem \ref{T2} just stops short of what is required. 
	
	An interesting question that we were not able \bblue{to} answer is whether or not the transformation $\Gamma$ from Theorem \ref{T3} is actually a bijection. If this were so then Theorem \ref{T1} would also yield local Lipschitz regularity of the gradient $DF$ outside a \bblue{locally} finite set of points in $\Omega$ for the differential inclusion $DF\in K$. This would be a very attractive result and would hint at the possibility of a regularity theory for differential inclusions into sets $S$ that do not have rank-$1$ connections but are not  elliptic.  \nl\nl
\em Acknowledgments. \rm The first author would like to thank Camillo DeLellis for pointing out that the Hilbert Schmidt norm of the matrix $M_h(x)$ (of the proof of Proposition 4.6. \cite{ADM}) tends towards $\frac{10}{36}$ as $h\rightarrow 0$. Roughly speaking $M_h(x)$ is (in the limit) analogous to $DF_u(x)$ from this paper and hence from this calculation it is clear that $F_u$ is a quasiregular mapping. Our desire to further investigate this observation was the starting point of this paper.   
The first author would also like to acknowledge the support of a Simons Foundation collobartive grant, award number 426900.

	\section{Sketch of the proof}
	
	As explained in the introduction, via the reduction to differential inclusions we get fractional Sobolev regularity $\na u \in W^{\bblue{\sigma},4}(\Omega')$ for all $\bblue{\sigma}\in (0,\frac{1}{3})$ \bblue{and all $\Omega'\subset\subset\Omega$}. In particular we have estimate  (\ref{eqss1}). The main thing we gain from this is the following estimate 
	(see Lemma \ref{l11}, (\ref{l114})) which is one of our key technical tools 
	\begin{equation}
	\label{l114.5}
	\int_{\Omega'}\lt|1-|\nabla u_{\ep}|^2\rt|\lt|u_{\ep,mn}\rt| \lt|f\rt| dx \leq C \lVert f \rVert_{L^{r}(\Omega')}\text{ for any }m,n\in \lt\{1,2\rt\}, r\geq 4.
	\end{equation}
	We will use (\ref{l114.5}) repeatedly.

	Our strategy will be to show that for 
	\bblue{\begin{equation}
	\label{s5.1}
	\Phi^{\xi}(z):=
	\begin{cases}
	|z|^2\xi &\text{ for } z\cdot\xi>0,\\
	0 &\text{ for } z\cdot\xi\leq 0,
	\end{cases}
	\end{equation}}
	we have 
	\begin{equation}
	\label{eqi922}
	\nabla\cdot\lt[\Phi^{\xi}\lt(\na u^{\perp}\rt)\rt]=0\text{ distributionally in }\Omega,\text{ for any }\xi\in \mathbb{S}^1\backslash \lt\{e_1,-e_1,e_2,-e_2\rt\}.
	\end{equation}  
	\bblue{Regularity} then follows by Theorem  \ref{JOP1} because any $\Phi^{\xi}\lt(\bblue{\na u(x)^{\perp}}\rt)=\xi \chi\lt(\bblue{x},\xi\rt)$ \bblue{for $|\na u(x)|=1$}, hence $\xi\cdot \na\chi\lt(\bblue{x},\xi\rt)=0$ distributionally in $\Omega$. This is a somewhat similar strategy to that of Ignat \cite{ignat} and DeLellis, Ignat \cite{DeI} except that in \cite{ignat}, \cite{DeI} it was shown that $\nabla\cdot\lt[\Phi\lt(\na u^{\perp}\rt)\rt]=0$ distributionally in $\Omega$ for all entropies $\Phi$, they then conclude (\ref{eqi45}) using (as explained in the introduction) the fact that $\xi \chi(\cdot,\xi)$ is the limit of a sequence of entropies. We will build toward establishing (\ref{eqi922}) in a couple of steps. \nl

	\em Step 1. \rm \em Harmonic entropies vanish: \rm In this step we  identify a class of entropies whose divergence vanishes when applied to $\na u^{\perp}$ as consequences of (\ref{eqi46}) holding. From Lemma 3 \cite{mul2} (see Lemma \ref{l13} \bblue{in this paper}) we know there is a one to one correspondence between entropies $\Phi$ and functions $\varphi\in C^{\infty}_c(\R^2)$ via the formula 
	\begin{equation}
	\label{eqi600}
	\Phi(z)=\varphi(z)z+\lt(\nabla\varphi(z)\cdot z^{\perp}\rt)z^{\perp}. 
	\end{equation}
	As we will sketch,  it will turn out that under the assumption of  (\ref{eqi46}),  if $\varphi$ is harmonic then $\nabla\cdot\lt[\Phi(\na u^{\perp})\rt]=0$. We will call entropies $\Phi$ that come from (\ref{eqi600}) via a harmonic $\varphi$,  \em harmonic entropies. \rm 
	
	%\footnote{Indeed its interesting that if we take the simplest harmonic polynomials $\varphi_1(z)=z_1^2-z_2^2$ and  $\varphi_2(z)=z_1 z_2$ then we recover entropies  $\wt{\Sigma}_{e_1 e_2}$ and $\wt{\Sigma}_{\ep_1 \ep_2}$ from the formula (\ref{eqi600}) }
	
	To see this we argue as follows. One of the key lemmas on entropies is Lemma 2 \cite{mul2} (see Lemma \bblue{\ref{l18} in this paper}), says that we can write 
	\begin{equation}
	\label{eqi801}
	\nabla\cdot\lt[\Phi(m)\rt]=\Psi(m)\cdot \na (1-\lt|m\rt|^2)\text{ for some }\Psi\in C_c^{\infty}(\R^2;\R^2).
	\end{equation}
	Now let $f_{\ep}:=f*\rho_{\ep}$ where $\rho_{\ep}(z)=\rho\lt(\frac{z}{\ep}\rt)\ep^{-2}$ and $\rho$ is the standard convolution kernel. Let $w=(w^1,w^2)=\na u^{\perp}$. \bblue{For $\Omega'\subset\subset\Omega$}, let $\zeta\in C^{\infty}_c(\Omega')$ be a test function, so integrating by parts we have 
	\begin{eqnarray*}
	%\label{eqi601}
	\int_{\Omega'} \nabla\cdot\lt[\Phi(w_{\ep})\rt] \zeta dx &\approx & -\int_{\Omega'} (1-\lt|w_{\ep}\rt|^2) \nabla\cdot\lt[\Psi(w_{\ep})\rt] \zeta dx\nn\\
	&=& -\int_{\Omega'} (1-\lt|w_{\ep}\rt|^2)\lt(\Psi_{1,1}(w_{\ep})w_{\ep,1}^1+\Psi_{1,2}(w_{\ep,1})w_{\ep,1}^2+\Psi_{2,1}(w_{\ep})w_{\ep,2}^1+\Psi_{2,2}(w_{\ep})w_{\ep,2}^2 \rt)\zeta dx.\nn\\
	\end{eqnarray*}
	The key point is that if $\Phi$ is a harmonic entropy then it is a calculation to see that  $\Psi_{1,2}=\Psi_{2,1}$. Now \bblue{we have} 
	\begin{equation}
	\label{eqi602}
	\nabla\cdot\lt[\wt{\Sigma}_{e_1 e_2}(\na u_{\ep}^{\perp})\rt]\overset{(\ref{l15.3})}{=}(u_{\ep,11}-u_{\ep,22})(1-\lt|\na u_{\ep}\rt|^2)=(w_{\ep,1}^2+w_{\ep,2}^1)(1-\lt|w_{\ep}\rt|^2)
	\end{equation}
	and 
	\begin{equation}
	\label{eqi603}
	\nabla\cdot\lt[\wt{\Sigma}_{\ep_1 \ep_2}(\na u_{\ep}^{\perp})\rt]\overset{(\ref{l15.4})}{=} 2 u_{\ep,12}(1-\lt|\na u_{\ep}\rt|^2)=-2 w^1_{\ep,1}(1-\lt|w_{\ep}\rt|^2)=2w^2_{\ep,2}(1-\lt|w_{\ep}\rt|^2).
	\end{equation}
	Proceeding formally and absorbing $\Psi_{1,1}(w_{\ep})$ into the test function $\zeta$ (strictly speaking we can not do this because $\Psi_{1,1}(w_{\ep})$ depends on $\ep$, however this can be overcome with estimate (\ref{l114.5})) we have that since $\nabla\cdot\lt[\wt{\Sigma}_{\ep_1 \ep_2}(\na u^{\perp})\rt]$ vanishes so
	$$
	\int_{\Omega'} (1-\lt|w_{\ep}\rt|^2)\Psi_{1,1}(w_{\ep})w^1_{\ep,1} \zeta dx\approx 0.
	$$ 
	In the same way  $\int_{\Omega'} (1-\lt|w_{\ep}\rt|^2)\Psi_{\bblue{2,2}}(w_{\ep})w^2_{\ep,2}\zeta dx\approx 0$ and, \bblue{since $\nabla\cdot\lt[\wt{\Sigma}_{e_1 e_2}(\na u^{\perp})\rt]=0$ and $\Psi_{1,2}=\Psi_{2,1}$},
	$$
	\int_{\Omega'} (1-\lt|w_{\ep}\rt|^2)\lt(\Psi_{1,2}(w_{\ep})w^2_{\ep,1}+\Psi_{2,1}(w_{\ep})w^1_{\ep,2} \rt)\zeta dx
	=\int_{\Omega'} (1-\lt|w_{\ep}\rt|^2)\Psi_{1,2}(w_{\ep})\bblue{\lt(w_{\ep,1}^2+w_{\ep,2}^1\rt)}\zeta dx\approx 0.
	$$
	Thus  $\nabla\cdot\lt[\Phi(\na u^{\perp})\rt]= 0$ for all harmonic entropies. \nl
	
	\em Step 2.  Estimate (\ref{eqi922}) holds: \rm As we can see the real issue of getting the divergences of entropies to vanish from hypothesis (\ref{eqi46}) is the term $ (1-\lt|w_{\ep}\rt|^2)\lt(\Psi_{1,2}(w_{\ep})w^2_{\ep,1}+\Psi_{2,1}(w_{\ep})w^1_{\ep,2} \rt)$. Given that we started with just two entropies $\wt{\Sigma}_{e_1 e_2}$ and $\wt{\Sigma}_{\ep_1 \ep_2}$ whose divergence vanishes (when applied to $\na u^{\perp}$) and end up with an entire class of entropies (what we call \em harmonic entropies\em) \bblue{whose} divergence vanishes, the natural way to proceed is to attempt to use our class of harmonic entropies  to further expand into a larger class of vanishing entropies. So what we need is a harmonic entropy to deal with terms of the form  $\Psi_{1,2}(w_{\ep})w^2_{\ep,1}+\Psi_{2,1}(w_{\ep})w^1_{\ep,2}$. It turns out there is a harmonic entropy that serves this purpose.

	Now notice that 
	\begin{eqnarray*}
	%\label{eqi677}
	&~&\int_{\Omega'} \lt(1-\lt|w_{\ep}\rt|^2\rt)\lt[\Psi_{1,2}(w_{\ep})w^2_{\ep,1}+\Psi_{2,1}(w_{\ep})w^1_{\ep,2}\rt] \zeta dx \nn\\
	&~&\qd\qd=\int_{\Omega'} \lt(1-\lt|w_{\ep}\rt|^2\rt)\frac{\lt(\Psi_{1,2}(w_{\ep})+\Psi_{2,1}(w_{\ep})\rt)}{2} \lt(w_{\ep,1}^2+w_{\ep,2}^1\rt) \zeta dx\nn\\
	&~&\qd\qd\qd\qd+\int_{\Omega'} \lt(1-\lt|w_{\ep}\rt|^2\rt)\frac{\lt(\Psi_{1,2}(w_{\ep})-\Psi_{2,1}(w_{\ep})\rt)}{2} \lt(w_{\ep,1}^2-w_{\ep,2}^1\rt) \zeta dx.
	\end{eqnarray*}
	The first term can be dealt with by absorbing $\frac{\lt(\Psi_{1,2}(w_{\ep})+\Psi_{2,1}(w_{\ep})\rt)}{2}$ into $\zeta$ as 		before then applying (\ref{eqi602}). So the term we have to deal with is the latter term. Now if $\varphi$ is related to $\Phi$ by 
	(\ref{eqi600}) it is a calculation to see that $\bblue{\Psi_{1,2}(z)-\Psi_{2,1}(z)\overset{(\ref{s4.2})}{=}\frac{1}{2} \na \lt( \Delta \varphi(z)\rt)\cdot z^{\perp}}=:		\psi(z)$. So 
	\begin{equation}
	\label{eqj2.88}
	\int_{\Omega'} \lt(1-\lt|w_{\ep}\rt|^2\rt)\frac{\lt(\Psi_{1,2}(w_{\ep})-\Psi_{2,1}(w_{\ep})\rt)}{2} \lt(w_{\ep,1}^2-w_{\ep,2}^1\rt) \zeta dx=			\frac{1}{2}\int_{\Omega'} \lt(1-\lt|w_{\ep}\rt|^2\rt)\psi\lt(w_{\ep}\rt)\lt(w_{\ep,1}^2-w_{\ep,2}^1\rt) \zeta dx,
	\end{equation}
	\xred{using} the fact $w_{\ep,1}^2-w_{\ep,2}^1=\bblue{\Delta u_{\ep}}$. Thus what we need is a harmonic entropy that includes the term $\bblue{\Delta u_{\ep}}$. Now taking $\varphi(z)=z_1^2-z_2^2$, via formula (\ref{eqi600}) we obtain entropy 
	$\Phi_0(z)=(z_1^3+3 z_1z_2^2,-3z_1^2z_2-z_2^3)$ and a short calculation gives 
	\begin{equation*}
	%\label{eqi688}
	\nabla\cdot[\Phi_0(w_{\ep})]=-6\lt(u_{\ep,1}u_{\ep,2}\bblue{\Delta u_{\ep}}+|\nabla u_{\ep}|^2u_{\ep,12}\rt).
	\end{equation*}
	Now \bblue{it is} a calculation (see (\ref{l912})) using (\ref{eqi801}) \bblue{to write}
	\begin{equation*}
	%\label{l912.6}
	\begin{split}
	&\int_{\Omega'}\lt(1-|\nabla u_{\ep}|^2\rt)\lt(u_{\ep,1}u_{\ep,2}\bblue{\Delta u_{\ep}}+|\nabla u_{\ep}|^2u_{\ep,12}\rt)\zeta dx\\
	=&-\frac{1}{12}\int_{\Omega'}  \na \cdot \lt[\Psi_0(w_{\ep}) \lt(1-\lt|w_{\ep}\rt|^2 \rt)^2  \rt] \zeta dx+\frac{1}{12}\int_{\Omega'}  \lt(1-\lt|w_{\ep}\rt|^2 \rt)^2 \na\cdot \lt[ \Psi_0(w_{\ep}) \rt] \zeta dx.
	\end{split}
	\end{equation*}
	The first term can be dealt with by integration by parts, and the second can be controlled via estimate (\ref{l114.5}). It follows that
	\begin{equation*}
	%\label{eqi900}
	\int_{\Omega'}\lt(1-|\nabla u_{\ep}|^2\rt)\lt(u_{\ep,1}u_{\ep,2}\bblue{\Delta u_{\ep}}+|\nabla u_{\ep}|^2u_{\ep,12}\rt)\zeta dx\rightarrow 0\text{ as }\ep\rightarrow 0. 
	\end{equation*}
	Now as $\lt|\na u\rt|=1$ a.e.\ we have 
	\begin{eqnarray}
	\label{eqj2.77}	
	\int_{\Omega'} \lt(1-\lt|\na u_{\ep}\rt|^2\rt)\lt|\na u_{\ep}\rt|^2 u_{\ep,12} \zeta dx &\approx & 
	\int_{\Omega'} \lt(1-\lt|\na u_{\ep}\rt|^2\rt) u_{\ep,12} \zeta dx\nn\\
	&\overset{(\ref{eqi603})}{=}&	\bblue{\frac{1}{2}\int_{\Omega'} \na\cdot\lt[\Sigma_{\ep_1 \ep_2} u_{\ep}\rt] \zeta dx}\rightarrow 0\qd\text{ as }\ep\rightarrow 0.\nn
	\end{eqnarray}
	So 
	\begin{equation}
	\label{eqj3}
	\int_{\Omega'} \lt(1-\lt|\na u_{\ep}\rt|^2\rt) u_{\ep,1} u_{\ep,2} \Delta u_{\ep}\bblue{\zeta} dx \rightarrow 0\qd\text{ as }\ep\rightarrow 0.
	\end{equation}
	Now from (\ref{eqj2.88}), \bblue{using $w_{\ep}=\na u_{\ep}^{\perp}$}, we can write 
	\begin{equation*}
	%\label{eqj2}
	\int_{\Omega'} \lt(1-\lt|w_{\ep}\rt|^2\rt)\frac{\lt(\Psi_{1,2}(w_{\ep})-\Psi_{2,1}(w_{\ep})\rt)}{2} \lt(w_{\ep,1}^2-w_{\ep,2}^1\rt) \zeta dx=\frac{1}{2}\int_{\Omega'} \lt(1-\lt|\na \bblue{u}_{\ep}\rt|^2\rt) u_{\ep,1} u_{\ep,2} \Delta u_{\ep}   \frac{\psi\lt(w_{\ep}\rt)\zeta}{ u_{\ep,1} u_{\ep,2}} dx.
	\end{equation*}
	\bblue{It turns out that} if 
	\begin{equation}
	\label{eqj4}
	\frac{\psi\lt(w_{\ep}\rt)}{ u_{\ep,1} u_{\ep,2}}\text{ remains uniformly bounded for small }\ep>0,
	\end{equation}
	then it can be absorbed (via estimate (\ref{l114.5})) into $\zeta$, \bblue{and as a result of (\ref{eqj3}) we have} 
	\begin{equation}
	\label{eqj6}
	\int_{\Omega'} \lt(1-\lt|w_{\ep}\rt|^2\rt)\frac{\lt(\Psi_{1,2}(w_{\ep})-\Psi_{2,1}(w_{\ep})\rt)}{2} \lt(w_{\ep,1}^2-w_{\ep,2}^1\rt) \zeta dx\rightarrow 0\qd\text{ as }\ep\rightarrow 0.
	\end{equation} 
	Hence the \bblue{estimate (\ref{eqj6}) holds} as long as (\ref{eqj4}) holds true. So we need to restrict ourselves to a class of entropies for which (\ref{eqj4}) is true. The key point is that for the sequence of entropies $\{\Phi^k\}$ that approximates $\Phi^{\xi}$ (for $\xi\in \mathbb{S}^1\backslash \lt\{e_1,-e_1,e_2,-e_2\rt\}$) we can guarantee that (\ref{eqj4}) holds true. Thus we can establish (\ref{eqi922}). \nl
	
	\em Sketch of proof completed. \rm The choice of coordinate system axis $\lt\{e_1,e_2\rt\}$ in (\ref{eqi922}) is completely arbitrary. We could have carried out the proof with the coordinate system axis $\lt\{\ep_1,\ep_2\rt\}$ and could then conclude (\ref{eqi922}) for any 
	$\bblue{\xi\in}\mathbb{S}^1\backslash \lt\{\ep_1,-\ep_1,\ep_2,-\ep_2\rt\}$. Thus (\ref{eqi922}) holds from any $\xi\in \mathbb{S}^1$ and therefore (\ref{eqi45}) holds true and 
	regularity follows by Theorem \ref{JOP1}.

%%%%%%%%%%%%%%%%%%%%%%%%%%%%%%%%%%%%%%%%%%%%%%%%%%%%%%
%
%
%%%%%%%%%%%%%%%%%%%%%%%%%%%%%%%%%%%%%%%%%%%%%%%%%%%%%%

\section{Background}
\label{background}

\bblue{In this section we provide} some background. Any two by two matrix can be uniquely decomposed into conformal and anticonformal parts as follows 
\begin{equation*}
\lt(\begin{matrix} a_{11} & a_{12} \\ a_{21} & a_{22}\end{matrix}\rt)=\frac{1}{2}
\lt(\begin{matrix} a_{11}+a_{22} & -(a_{21}-a_{12}) \\ a_{21}-a_{12} & a_{11}+a_{22}\end{matrix}\rt)
+\frac{1}{2}
\lt(\begin{matrix} a_{11}-a_{22} & a_{21}+a_{12} \\ a_{21}+a_{12} & -(a_{11}-a_{22})\end{matrix}\rt).
\end{equation*}
So for a matrix \xblue{$A=\lt(\begin{matrix} a_{11} & a_{12} \\ a_{21} & a_{22}\end{matrix}\rt)$}, define 
\begin{equation}\label{eqzzz8}
\xblue{\lt[A\rt]_c:=\frac{1}{2}
\lt(\begin{matrix} a_{11}+a_{22} & -(a_{21}-a_{12}) \\ a_{21}-a_{12} & a_{11}+a_{22}\end{matrix}\rt)}
\text{ and }\lt[A\rt]_a:=\frac{1}{2}
\lt(\begin{matrix} a_{11}-a_{22} & a_{21}+a_{12} \\ a_{21}+a_{12} & -(a_{11}-a_{22})\end{matrix}\rt).
\end{equation}
Its easy to see that 
\begin{equation}
\label{eqzzz9}
\det\lt(A\rt)=\det(\lt[A\rt]_c)+\det(\lt[A\rt]_a).
\end{equation}

Given $w:\Omega\rightarrow \R^2$ \xblue{such that} $w(\xblue{x_1,x_2})=(u(\xblue{x_1,x_2}),v(\xblue{x_1,x_2}))$, for 
$z=x_1+ix_2$, let $\wii(z)=u(\xblue{x_1,x_2})+i v(\xblue{x_1,x_2})$. Note that
$\frac{\partial \wii}{\partial \overline{z}}(z)=\frac{1}{2}(\frac{\partial}{\partial x_1}+i\frac{\partial}{\partial x_2})\wii=\frac{1}{2}(u_{\xblue{,1}}-v_{\xblue{,2}})+\frac{i}{2}(v_{\xblue{,1}}+u_{\xblue{,2}})$ and $\frac{\partial \wii}{\partial z}(z)=\frac{1}{2}(\frac{\partial}{\partial x_1}-i\frac{\partial}{\partial x_2})\wii=\frac{1}{2}(u_{\xblue{,1}}+v_{\xblue{,2}})+\frac{i}{2}(v_{\xblue{,1}}-u_{\xblue{,2}})$. Now 
identifying complex numbers with conformal matrices in the standard way 
\begin{equation}
\label{eqdf3}
\lt[x_1+ix_2\rt]_M=\lt(\begin{matrix} x_1 & -x_2 \\ x_2 & x_1\end{matrix}\rt),
\end{equation}
we have 
that 
\begin{equation}
\label{eqzzz2}
\lt[D w(x)\rt]_a=\blue{\lt[\frac{\partial \wii}{\partial \overline{z}}(z)\rt]_M
	\lt(\begin{matrix} 1 & 0 \\ 0 & -1\end{matrix}\rt)}\text{ and }\lt[D w(x)\rt]_c=\blue{\lt[\frac{\partial \wii}{\partial z}(z)\rt]_M}.
\end{equation}

\section{ Proof of Theorem \ref{T3}}

\begin{a1}\label{L1}
	Let \xblue{$\Omega$ and $\wt{\Omega}$ be as in Theorem \ref{T3}. Define}
	\begin{equation}
	\label{eqy29}
	K:=\lt\{\lt(\begin{array}{cc} \frac{2}{3}\sin^3(\theta) &  \frac{2}{3}\cos^3(\theta)  \\ -\cos(\theta)\lt(1-\frac{2}{3}\cos^2(\theta)\rt)  & \sin(\theta)\lt(1-\frac{2}{3}\sin^2(\theta)\rt)\end{array}\rt):\theta\in \lt[0,2\pi\rt)  \rt\}.
	\end{equation}
	\xblue{Let a map $F=(F_1,F_2) \in W^{1,\infty}(\Omega;\Rb^2)$ and a function $v\in W^{1,\infty}(\wt\Omega;\mathbb{C})$ be related by $v(\xblue{x_1+ix_2})=F_1(\xblue{x_1,x_2})+i F_2(\xblue{x_1,x_2})$. Then $DF\in K$ at $x\in\Omega$ if and only if $v$ satisfies the following non-linear Beltrami equation and constraint at $z=x_1+ix_2\in\wt\Omega$:
	\begin{equation}
	\label{eqy20}
	\frac{\partial v}{\partial \bar{z}}(\xblue{z})=\frac{4}{3}\lt(\frac{\partial v}{\partial z}(\xblue{z}) \rt)^3,    \qd\lt|\frac{\partial v}{\partial z}(\xblue{z})\rt|=\frac{1}{2}.
	\end{equation}}
\end{a1}

\begin{proof}[Proof of Lemma \ref{L1}]
	First \xblue{assume $x\in\Omega$ is such that $DF(x)\in K$. We show that $v$ satisfies (\ref{eqy20}) at $z=x_1+ix_2$. } Note that since $DF(x)\in K$, there exists $\theta\in\lt[0,2\pi\rt)$ such that 
	\begin{equation}
	\label{eqz19}
	DF(x)= \lt(\begin{array}{cc}   \frac{2}{3}\sin^3(\theta)    & \frac{2}{3}\cos^3(\theta) \\ -\cos(\theta)\lt(1-\frac{2}{3}\cos^2(\theta)\rt)  &    \sin(\theta)\lt(1-\frac{2}{3}\sin^2(\theta)\rt)        \end{array}\rt).
	\end{equation}
	As described in Section \ref{background}, for any matrix $A$, we decompose $A = [A]_c + [A]_a$, where $[A]_c$ and $[A]_a$ are the conformal and anticonformal parts of $A$, respectively. Using (\ref{eqzzz8}) and (\ref{eqz19}) we have
	\begin{equation}
	\label{eqz20}
	\lt[DF(x)\rt]_c=\frac{1}{2}\lt(\begin{array}{cc} \sin(\theta)  & \cos(\theta) \\ -\cos(\theta) & \sin(\theta)\end{array}\rt). 
	\end{equation}
	Now recalling the trig identities 
	\begin{equation}
	\label{eqz20.5}
	\sin\lt(3\theta\rt)=-4\sin^3\lt(\theta\rt)+3\sin\lt(\theta\rt),\; \cos\lt(3\theta\rt)=4\cos^3\lt(\theta\rt)-3\cos\lt(\theta\rt).
	\end{equation}
	Note that
	\begin{eqnarray}
	\label{eqdf2}
	\lt[DF(x)\rt]_a&=&\frac{1}{2}\lt(\begin{array}{cc} \frac{4}{3}\sin^3(\theta)-\sin(\theta)  & \frac{4}{3}\cos^3(\theta)- \cos(\theta) \\  \frac{4}{3}\cos^3(\theta)-\cos(\theta) &  -\lt(   \frac{4}{3}\sin^3(\theta)-\sin(\theta) \rt)   \end{array}\rt)\nn\\
	&\overset{(\ref{eqz20.5})}{=}&\frac{1}{2}\lt(\begin{array}{cc} -\frac{1}{3}\sin(3\theta) & \frac{1}{3}\cos(3\theta) \\ \frac{1}{3}\cos(3\theta) & \frac{1}{3}\sin(3\theta) \end{array}\rt)\nn\\
	&=&\frac{1}{6}\lt(\begin{array}{cc} -\sin(3\theta) & \cos(3\theta) \\ \cos(3\theta) &\sin(3\theta) \end{array}\rt).
	\end{eqnarray}
	Recall that $v(x_1+ix_2)=F_1(x_1,x_2)+iF_2(x_1,x_2)$. It follows from (\ref{eqzzz2}) that
	\begin{equation}
	\label{eqdf1}
	\lt[\frac{\partial v}{\partial z}(z)\rt]_M\overset{(\ref{eqzzz2})}{=}\lt[ D F(x)\rt]_c \qd\text{ and } \qd\lt[\frac{\partial v}{\partial \bar{z}}(z)\rt]_M\overset{(\ref{eqzzz2})}{=}\lt[ D F(x)\rt]_a\lt(\begin{array}{cc}   1  & 0 \\
	0 &  -1 \end{array}\rt). 
	\end{equation}
	Thus 
	\begin{equation}\label{eqz1}
	\frac{\partial v}{\partial z}(z)\overset{(\ref{eqdf3}),(\ref{eqz20}),(\ref{eqdf1})}{=}\frac{1}{2}\lt(\sin(\theta)-i\cos(\theta)\rt),
	\end{equation}
	\begin{eqnarray}
	\label{eqdf4}
	\lt[ D F(x)\rt]_a\lt(\begin{array}{cc}   1  & 0 \\
	0 &  -1 \end{array}\rt)&\overset{(\ref{eqdf2})}{=}&\frac{1}{6} \lt(\begin{array}{cc}   -\sin(3\theta)  & -\cos(3\theta) \\
	\cos(3\theta) &  -\sin(3\theta)  \end{array}\rt).
	\end{eqnarray}
	Therefore it follows that
	\begin{eqnarray}\label{eqs21}
	\frac{\partial v}{\partial \bar{z}}(z)&\overset{(\ref{eqdf1}),(\ref{eqdf4}),(\ref{eqdf3})}{=}&-\frac{1}{6}\lt(\sin(3\theta)-i\cos(3\theta)\rt)\nn\\
	&=&\frac{1}{6}\lt(\sin(\theta)-i\cos(\theta)\rt)^3\nn\\
	&\overset{(\ref{eqz1})}{=}&\frac{1}{6}\lt(2 \frac{\partial v}{\partial z}(z) \rt)^3\nn\\
	&=&\frac{4}{3}\lt(\frac{\partial v}{\partial z}(z) \rt)^3.
	\end{eqnarray}
	\xblue{We obtain from (\ref{eqs21}) and (\ref{eqz1}) that $v$ satisfies the constrained non-linear Beltrami equation (\ref{eqy20}) at $z\in\wt\Omega$.}
	
	Conversely, suppose the function $v\in W^{1,\infty}(\wt\Omega;\Cb)$ satisfies (\ref{eqy20}) at $z=x_1+ix_2$. Recall that $F(x_1,x_2)=\lt(\mathrm{Re}(v(x_1+ix_2)),\mathrm{Im}(v(x_1+ix_2))\rt)$. We will show that $DF(x) \in K$. Indeed, we have
	\begin{equation}
	\label{eqy3}
	\frac{\partial v}{\partial z}=\frac{1}{2}\bigg[\lt(F_{1,1}+F_{2,2}\rt)+i\lt(F_{2,1}-F_{1,2}\rt)\bigg]
	\end{equation}
	and
	\begin{equation}
	\frac{\partial v}{\partial \bar{z}}=\frac{1}{2}\bigg[\lt(F_{1,1}-F_{2,2}\rt)+i\lt(F_{2,1}+F_{1,2}\rt)\bigg].
	\end{equation}
	Since $\lt|\frac{\partial v}{\partial z}(z) \rt|=\frac{1}{2}$, there exists $\theta\in[0,2\pi)$ such that
	\begin{equation}
	\frac{\partial v}{\partial z} = \frac{1}{2}\lt(\cos(\theta)+i\sin(\theta)\rt).
	\end{equation}
	\xblue{Now since $v$ satisfies (\ref{eqy20}) at $z$, we have}
	\begin{eqnarray}
	\label{eqy4}
	\frac{\partial v}{\partial \bar{z}}(z)&=&\frac{4}{3}\lt(\frac{1}{2}\lt(\cos(\theta)+i\sin(\theta)\rt)\rt)^3\nn\\
	&=&\frac{1}{6}\bigg(\cos(3\theta)+i\sin(3\theta)\bigg).
	\end{eqnarray}
	Now we obtain from (\ref{eqy3})-(\ref{eqy4}) that
		\begin{equation}
		\label{eqds15}
		\begin{split}
		F_{1,1}+F_{2,2}&=\cos(\theta),\\
		F_{2,1}-F_{1,2}&=\sin(\theta),\\
		F_{1,1}-F_{2,2}&=\frac{1}{3}\cos(3\theta),\\
		F_{2,1}+F_{1,2}&=\frac{1}{3}\sin(3\theta).
		\end{split}
		\end{equation}
		So solving (\ref{eqds15}) for $F_{1,1}, F_{1,2}, F_{2,1}, F_{2,2}$,  we obtain
	\begin{equation*}
	\begin{split}
	F_{1,1}&=\frac{1}{2}\cos(\theta)+\frac{1}{6}\cos(3\theta)\overset{(\ref{eqz20.5})}{=}\frac{2}{3}\cos^3(\theta),\\
	F_{1,2}&=\blue{\frac{1}{6}\sin(3\theta)-\frac{1}{2}\sin(\theta)} \overset{(\ref{eqz20.5})}{=} -\frac{2}{3}\sin^3(\theta),\\
	F_{2,1}&=\frac{1}{2}\sin(\theta)+\frac{1}{6}\sin(3\theta) \overset{(\ref{eqz20.5})}{=} \sin(\theta)-\frac{2}{3}\sin^3(\theta),\\
	F_{2,2}&=\frac{1}{2}\cos(\theta)-\frac{1}{6}\cos(3\theta)\overset{(\ref{eqz20.5})}{=} \cos(\theta)-\frac{2}{3}\cos^3(\theta).
	\end{split}
	\end{equation*}
	Now letting $\tilde{\theta}=\frac{\pi}{2}+\theta$, we have $\cos(\tilde{\theta})=-\sin(\theta)$ and $\sin(\tilde{\theta})=\cos(\theta)$. One can check immediately that $DF \in K$ at $x=(x_1,x_2)$ with the phase function $\wt{\theta}$.
\end{proof}

%%%%%%%%%%%%%%%%%%%%%%%%%%%%%%%%%%%%%%%%%%%%%%%%%%%
%
%
%%%%%%%%%%%%%%%%%%%%%%%%%%%%%%%%%%%%%%%%%%%%%%%%%%%%%%%%

\subsection{Proof of Theorem \ref{T3} completed} Firstly given $u\in E(\Omega)$ we can define $F_u:\Omega\rightarrow \R^2$ by $F_u(0,0)=(0,0)$ and 
\begin{equation}
\label{eqds20}
DF_u=\lt(\begin{array}{cc} u_{,2}\lt(1-u_{,1}^2-\frac{u_{,2}^2}{3}\rt)  & u_{,1}\lt(1-u_{,2}^2-\frac{u_{,1}^2}{3}\rt)  \\ -u_{,1}\lt(1-\frac{2 u_{,1}^2}{3}\rt)  &  u_{,2}\lt(1-\frac{2 u_{,2}^2}{3}\rt) \end{array}\rt).
\end{equation}
\green{The existence of $F_u$ over bounded simply-connected Lipschitz domains in the classical $L^2$ framework can be found in \cite{gr}. We provide a proof of the existence of $F_u$ over bounded simply-connected domains in Lemma \ref{curl} in the Appendix. Such results might be well-known to experts, but we were not able to find a reference. Therefore we include a proof for the convenience of the readers.} \xblue{Since $|\nabla u|=1$ a.e. in $\Omega$, it is clear that $F_u\in W^{1,\infty}(\Omega;\Rb^2)$} is a mapping that satisfies 
\begin{equation*}
%\label{eq18}
DF_u \in \lt\{\lt(\begin{array}{cc} \sin(\theta)(1-\cos^2(\theta)-\frac{\sin^2(\theta)}{3})  & \cos(\theta)(1-\sin^2(\theta)-\frac{\cos^2(\theta)}{3})  \\ -\cos(\theta)\lt(1-\frac{2}{3}\cos^2(\theta)\rt)  & \sin(\theta)\lt(1-\frac{2}{3}\sin^2(\theta)\rt)\end{array}\rt)  :\theta\in \lt[0,2\pi\rt)\rt\}\overset{(\ref{eqy29})}{=}K \xblue{\text{ a.e. in } \Omega.} 
\end{equation*} 
Thus applying Lemma \ref{L1} we have that $v_u(x_1+ix_2):=F_u^1(x_1,x_2)+iF_u^2(x_1,x_2)$ satisfies the non-linear Beltrami system (\ref{eqs6.5}). 
So defining 
\begin{equation}\label{eqz2}
\Gamma(u):=v_u,
\end{equation} 
we have that $\Gamma$ forms a transformation of $\lt[E(\Omega)/\xblue{\Rb}\rt]$ into $B(\wt{\Omega})$. Now we show that $\Gamma$ is injective. \xblue{Given $u,w\in \lt[E(\Omega)/\xblue{\Rb}\rt]$ such that $\Gamma(u)=\Gamma(w)$, we have $DF_u=DF_w$. Note that for all $x\in\Omega$ such that $|\nabla u(x)|=1$, we deduce from (\ref{eqds20}) that
\begin{equation}\label{eqz59}
u_{,1}=F_{u,2}^1-F_{u,1}^2 \qd\text{ and }\qd u_{,2}=F_{u,1}^1+F_{u,2}^2.
\end{equation}
The same relations hold for $\nabla w$. This implies $\na u=\na w$ a.e. in $\Omega$ and hence $u=w$ in $\lt[E(\Omega)/\xblue{\Rb}\rt]$. Thus we have shown that $\Gamma$ is injective. } 

Now for the second part of the theorem, given a function $v\in B(\wt{\Omega})\cap W^{2,1}(\wt{\Omega})$ we need to show that there exists some $u\in \lt[E(\Omega)/\Rb\rt]\cap W^{2,1}(\wt{\Omega})$ such that $\Gamma(u)=v$. Let us define 
\begin{equation*}
%\label{eq18.7}
\bblue{F}(x_1,x_2)=\lt(\mathrm{Re}(v(x_1+ix_2)),\mathrm{Im}(v(x_1+ix_2))\rt). 
\end{equation*} 
By Lemma \ref{L1} we have  
\begin{equation*}
%\label{eq19}
D\bblue{F}\in K\; a.e. \text{ in }\Omega. 
\end{equation*}
We have that $D\bblue{F}\in W^{1,1}(\Omega)$ and there exists $\theta(x):\Omega\rightarrow\lt[0,2\pi\rt)$ such that 
\begin{equation}
\label{eq20}
D\bblue{F}(x)=\lt(\begin{array}{cc} \sin(\theta(x))\lt(1-\cos^2(\theta(x))-\frac{\sin^2(\theta(x))}{3}\rt)  & \cos(\theta(x))\lt(1-\sin^2(\theta(x))-\frac{\cos^2(\theta(x))}{3}\rt)  \\ -\cos(\theta(x))\lt(1-\frac{2}{3}\cos^2(\theta(x))\rt)  & \sin(\theta(x))\lt(1-\frac{2}{3}\sin^2(\theta(x))\rt)\end{array}\rt)
\end{equation}
for a.e. $x\in \Omega$. \xblue{Similar to (\ref{eqz59}), we deduce from (\ref{eq20}) that
\begin{equation*}
	\cos\lt(\theta(x)\rt)=\bblue{F}_{1,2}(x)-\bblue{F}_{2,1}(x) \qd\text{ and }\qd \sin(\theta(x))=\bblue{F}_{1,1}(x)+\bblue{F}_{2,2}(x) \qd\text{a.e. in } \Omega.
\end{equation*}}
	Hence $\alpha(x):=\cos\lt(\theta(x)\rt)$ and $\beta(x):=\sin\lt(\theta(x)\rt)$ are such that $\alpha,\beta\in W^{1,1}(\Omega)$. Now we have, \xblue{for a.e. $x\in\Omega$,}
	\begin{eqnarray}
	\label{eqb18}
	&~&0=\xblue{\mathrm{curl}\lt(\nabla \bblue{F}_1\rt)=\mathrm{curl}\lt(\beta(x)\lt(1-\alpha(x)^2-\frac{\beta(x)^2}{3}\rt),\alpha(x)\lt(1-\beta(x)^2-\frac{\alpha(x)^2}{3}\rt)\rt)}\nn\\
	&~&\qd =
	\lt(1-\alpha(x)^2-\beta(x)^2\rt)\lt(\alpha_{,1}(x)-\blue{\beta_{,2}(x)}\rt)
	+2\alpha(x)\beta(x)\lt(\alpha_{,2}(x)-\beta_{,1}(x)\rt)\nn\\
	&~&\qd =2\alpha(x)\beta(x)\lt(\alpha_{,2}(x)-\beta_{,1}(x)\rt),
	\end{eqnarray}
    and
	\begin{eqnarray}
	\label{eqb60}
	&~&0=\xblue{\mathrm{curl}\lt(\nabla \bblue{F}_2\rt)=\mathrm{curl}\lt(-\alpha(x)\lt(1-\frac{2}{3}\alpha(x)^2\rt), \beta(x)\lt(1-\frac{2}{3}\beta(x)^2\rt)\rt)}\nn\\
	&~&\qd =\beta_{,1}(x)\lt(1-2\beta(x)^2\rt)+\alpha_{,2}(x)\lt(1-2\alpha(x)^2\rt)\nn\\
	&~&\qd=\beta_{,1}(x)\lt(1-2\beta(x)^2\rt)+\beta_{,1}(x)\lt(1-2\alpha(x)^2\rt)+\lt(\alpha_{,2}(x)-\beta_{,1}(x)\rt)\lt(1-2\alpha(x)^2\rt)\nn\\
	&~&\qd=2 \beta_{,1}(x)\lt(1-\beta(x)^2-\alpha(x)^2\rt)+\lt(\alpha_{,2}(x)-\beta_{,1}(x) \rt)\lt(\beta(x)^2-\alpha(x)^2\rt)\nn\\
	&~&\qd=\lt(\alpha_{,2}(x)-\beta_{,1}(x) \rt)\lt(\beta(x)^2-\alpha(x)^2\rt).
	\end{eqnarray}
	Taking the squares of (\ref{eqb18}) and (\ref{eqb60}) and adding, and using the fact that $\alpha(x)^2+\beta(x)^2=1$, we have
	\begin{equation*}
		\begin{split}
			0 &= \lt[\lt(\alpha(x)^2-\beta(x)^2\rt)^2 + 4\alpha(x)^2\beta(x)^2\rt] \lt|\mathrm{curl}(\alpha(x),\beta(x))\rt|^2\\
			&=\lt(\alpha(x)^2+\beta(x)^2\rt)^2 \lt|\mathrm{curl}(\alpha(x),\beta(x))\rt|^2 = \lt|\mathrm{curl}(\alpha(x),\beta(x))\rt|^2\text{ for a.e. } x\in \Omega.
		\end{split}
	\end{equation*}
	Therefore, we have
	\begin{equation*}
	\mathrm{curl}\lt(\alpha(x),\beta(x)\rt) = 0  \qd\text{for a.e. } x\in\Omega.
	\end{equation*}
	\xblue{Since $(\alpha(x),\beta(x))\in L^{\green{\infty}}(\Omega;\Rb^2)$, by \green{Lemma \ref{curl} in the Appendix}, there exists $u\in H^1(\Omega)$ such that $\nabla u=(\alpha,\beta)=(\cos(\theta),\sin(\theta))$. Since $\alpha(x)^2+\beta(x)^2=1$, it is clear that $u$ also belongs to $W^{1,\infty}(\Omega)$. This along with (\ref{eq20}) and the fact that $\alpha,\beta\in W^{1,1}(\Omega)$ implies that $u\in \lt[E(\Omega)/\Rb\rt]\cap W^{2,1}(\Omega)$. Now looking at (\ref{eq20}) and the definition of $\Gamma$ in (\ref{eqz2}), it is clear that $\Gamma(u)=v$. Hence, this completes the proof of the bijective part of the theorem. }

%%%%%%%%%%%%%%%%%%%%%%%%%%%%%%%%%%%%%%%%%%%%%%%%%%%%%%%%%%%
%
%
%%%%%%%%%%%%%%%%%%%%%%%%%%%%%%%

%
\section{Proof of Theorem \ref{T2}} Define $\blue{F(x_1,x_2)}=\lt(\mathrm{Re}\lt(v(x_1+ix_2)\rt),\mathrm{Im}\lt(v(x_1+ix_2)\rt)\rt)$. By Lemma \ref{L1} the function $F$ satisfies the differential inclusion $DF\in K$ a.e.\ in $\Omega$, where $K$ is the subset of all two by two matrices defined by (\ref{eqy29}). Let $\blue{\Omega=\lt\{(x_1,x_2):x_1+ix_2\in \wt\Omega\rt\}}$. Define 
\begin{equation*}
%\label{eqx100}
M(\theta):= \lt(\begin{array}{cc}   \frac{2}{3}\sin^3(\theta)    & \frac{2}{3}\cos^3(\theta) \\ -\cos(\theta)\lt(1-\frac{2}{3}\cos^2(\theta)\rt)  &    \sin(\theta)\lt(1-\frac{2}{3}\sin^2(\theta)\rt)        \end{array}\rt).
\end{equation*}
By Lemma \ref{L1}, there exists $\psi:\Omega\rightarrow\lt[0,2\pi\rt)$ such that
\begin{equation*}
%\label{eqy61}
DF(x)=M\lt(\psi(x)\rt) \text{ for a.e. } x\in \Omega.
\end{equation*}
\bblue{Given $\Omega'\subset\subset\Omega$, denote $\gamma:=\mathrm{dist}(\Omega',\partial\Omega)>0$. } Let $h\in B_{\bblue{\gamma}}(0)$ and define 
\begin{equation}
\label{eqy62}
\alpha_h(x)=\psi(x+h)-\psi(x)\text{ for } \blue{x\in \Omega'}.
\end{equation}
\xblue{First, we prove the following lemma.}

\xblue{\begin{a1}\label{l7}
For all $x\in\Omega'$ and $h\in B_{\bblue{\gamma}}(0)$ such that $DF(x),DF(x+h)\in K$, we have that
\begin{equation}\label{eqzz41}
		\det\lt(DF(x+h)-DF(x)\rt)> c_0\lt|DF(x+h)-DF(x)\rt|^4,
		\end{equation}
		where the constant $c_0$ is independent of $x$ and $h$.
\end{a1}}

\begin{proof}
\xblue{Given $x\in\Omega'$ and $h\in B_{\bblue{\gamma}}(0)$ such that $DF(x),DF(x+h)\in K$, we will show the estimate (\ref{eqzz41}) in several steps.}

\em Step 1. \rm We have
	\begin{equation}
	\label{eqs100}
	\det\lt(DF(x+h)-DF(x)\rt)=\frac{4}{9}-\frac{2}{3}\cos(\alpha_h(x))+\frac{2}{9}\cos^3(\alpha_h(x))=\frac{\alpha_h^4}{6}+o(\alpha_h^4).
	\end{equation}

\em Proof of Step 1. \rm We know 
\begin{eqnarray}
\label{eqs52}
DF(x)&=&\lt[DF(x)\rt]_c+\lt[DF(x)\rt]_a\nn\\
&\overset{(\ref{eqz20}),(\ref{eqdf2})}{=}&\frac{1}{2}\lt(\begin{array}{cc}  \sin(\psi(x))  & \cos(\psi(x)) \\
-\cos(\psi(x)) &  \sin(\psi(x)) \end{array}\rt)+\frac{1}{6} \lt(\begin{array}{cc}   -\sin(3\psi(x))  & \cos(3\psi(x)) \\
\cos(3\psi(x)) &  \sin(3\psi(x))  \end{array}\rt).
\end{eqnarray}
It follows that 
\begin{eqnarray}
\label{eqs2}
&~&DF(x+h)-DF(x)\nn\\
&~&\qd\qd \overset{(\ref{eqs52})}{=}\frac{1}{2}\lt(\begin{array}{cc}  \sin(\psi(x+h))-\sin(\psi(x))  & \cos(\psi(x+h))-\cos(\psi(x)) \\
-\cos(\psi(x+h))+\cos(\psi(x)) &  \sin(\psi(x+h))-\sin(\psi(x))  \end{array}\rt)\nn\\
&~&\qd\qd\qd\qd+\frac{1}{6} \lt(\begin{array}{cc}   -\sin(3\psi(x+h))+\sin(3\psi(x))  & \cos(3\psi(x+h))-\cos(3\psi(x)) \\
\cos(3\psi(x+h))-\cos(3\psi(x)) &  \sin(3\psi(x+h))-\sin(3\psi(x))  \end{array}\rt).
\end{eqnarray}
So \blue{using (\ref{eqzzz9}) we have}
\begin{eqnarray}
\label{eqs3}
&~&\det\lt(DF(x+h)-DF(x)\rt)\nn\\
&~&\qd\qd\qd=\frac{1}{4}\lt(\lt( \sin(\psi(x+h))-\sin(\psi(x))\rt)^2+\lt(\cos(\psi(x+h))-\cos(\psi(x))\rt)^2\rt)\nn\\
&~&\qd\qd\qd\qd-\frac{1}{36}\lt(\lt( \sin(3\psi(x+h))-\sin(3\psi(x))\rt)^2+\lt(\cos(3\psi(x+h))-\cos(3\psi(x))\rt)^2\rt)\nn\\
&~&\qd\qd\qd=\frac{1}{4}\lt(2-2 \sin(\psi(x+h))\sin(\psi(x))-2\cos(\psi(x+h))\cos(\psi(x))\rt)\nn\\
&~&\qd\qd\qd\qd-\frac{1}{36}\lt(2-2 \sin(3\psi(x+h))\sin(3\psi(x))-2\cos(3\psi(x+h))\cos(3\psi(x))\rt).
\end{eqnarray}
Recall that $\alpha_h(x)$ is defined by (\ref{eqy62}), so $\psi(x+h)=\psi(x)+\alpha_h(x)$. Now 
\begin{equation}
\label{eqs4}
\sin\lt(\psi(x+h)\rt)=\sin(\psi(x))\cos(\alpha_h(x))+\cos(\psi(x))\sin(\alpha_h(x))
\end{equation}
and
\begin{equation}
\label{eqs5}
\cos\lt(\psi(x+h)\rt)=\cos(\psi(x))\cos(\alpha_h(x))-\sin(\psi(x))\sin(\alpha_h(x)).
\end{equation}
Thus 
\begin{eqnarray}
\label{eqs6}
&~&2-2 \sin(\psi(x+h))\sin(\psi(x))-2\cos(\psi(x+h))\cos(\psi(x))\nn\\
&~&\qd\qd\overset{(\ref{eqs4}),(\ref{eqs5})}{=}2-2\lt(\sin(\psi(x))\cos(\alpha_h(x))+\cos(\psi(x))\sin(\alpha_h(x))\rt)\sin(\psi(x))\nn\\
&~&\qd\qd\qd\qd-2\lt(\cos(\psi(x))\cos(\alpha_h(x))-\sin(\psi(x))\sin(\alpha_h(x))\rt)\cos(\psi(x))\nn\\
&~&\qd\qd =2\lt(1-\cos(\alpha_h(x))\rt).
\end{eqnarray}
Note that $3\psi(x+h)-3\psi(x)=3\alpha_h(x)$, so $3\psi(x+h)=3\psi(x)+3\alpha_h(x)$. Thus 
\begin{eqnarray}
\label{eqs7}
&~&2-2\sin(3\psi(x+h))\sin(3\psi(x))-2\cos(3\psi(x+h))\cos(3\psi(x))\nn\\
&~&\qd\qd=2-2\lt(\sin(3\psi(x))\cos(3\alpha_h(x))+\cos(3\psi(x))\sin(3\alpha_h(x))\rt)\sin(3\psi(x))\nn\\
&~&\qd\qd\qd-2\lt(\cos(3\psi(x))\cos(3\alpha_h(x))-\sin(3\psi(x))\sin(3\alpha_h(x))\rt)\cos(3\psi(x))\nn\\
&~&\qd\qd=2\lt(1-\cos(3\alpha_h(x))\rt).
\end{eqnarray}
Thus putting (\ref{eqs6}) and (\ref{eqs7}) together with (\ref{eqs3}) we have that 
\begin{equation*}
%\label{eqs8}
\det\lt(DF(x+h)-DF(x)\rt)=\frac{1}{2}\lt(1-\cos(\alpha_h(x))\rt)-\frac{1}{18}\lt(1-\cos(3\alpha_h(x))\rt). 
\end{equation*}
Now $\blue{\cos\lt(3\alpha_h(x)\rt)\overset{(\ref{eqz20.5})}{=}4\cos^3\lt(\alpha_h(x)\rt)-3\cos\lt(\alpha_h(x)\rt)}$. So 
\begin{eqnarray*}
	\det\lt(DF(x+h)-DF(x)\rt)&=&\frac{1}{2}\lt(1-\cos(\alpha_h(x))\rt)-\frac{1}{18}\lt(1-4 \cos^3(\alpha_h(x))+3 \cos(\alpha_h(x))\rt)\nn\\
	&=&\frac{4}{9}-\frac{2}{3}\cos(\alpha_h(x))+\frac{2}{9}\cos^3(\alpha_h(x)).
\end{eqnarray*}
Now since $\cos(\alpha_h(x))= 1-\frac{\alpha_h^2}{2}+\frac{\alpha_h^4}{\blue{24}}+o(\alpha_h^4)$, we have
\begin{eqnarray*}
	\frac{4}{9}-\frac{2}{3}\cos(\alpha_h(x))+\frac{2}{9}\cos^3(\alpha_h(x))&=& \frac{4}{9}-\frac{2}{3}\lt(1-\frac{\alpha_h^2}{2}+\frac{\alpha_h^4}{\blue{24}}\rt)+\frac{2}{9}\lt(1-\frac{\alpha_h^2}{2}+\frac{\alpha_h^4}{\blue{24}}\rt)^3+o(\alpha_h^4)\nn\\
	&=&-\frac{2}{9}+\frac{\alpha_h^2}{3}-\frac{\alpha_h^4}{\blue{36}}+\frac{2}{9}\lt(1-\frac{3}{2}\alpha_h^2+\blue{\frac{7}{8}}\alpha_h^4\rt)+o(\alpha_h^4)\nn\\
	&=&\frac{\alpha_h^4}{6}+o(\alpha_h^4)
\end{eqnarray*}
\xblue{for $\alpha_h>0$ sufficiently small.}\nl

\blue{\em Step 2. \rm We have
	\begin{equation}\label{eqs10}
	\lt|DF(x+h)-DF(x)\rt|^2 = \frac{10}{9}-\frac{2}{3}\cos(\alpha_h(x))-\frac{4}{9}\cos^3(\alpha_h(x)) = \alpha_h^2+o(\alpha_h^2).
	\end{equation}}

\blue{\em Proof of Step 2.} \rm Now looking at \eqref{eqs2}, it is clear that the two matrices in the decomposition are orthogonal when they are identified as vectors in $\mathbb{R}^4$. \red{Therefore, using similar calculations as \blue{in Step 1}, we have
	\begin{equation*}
		\begin{split}
			&\lt|DF(x+h)-DF(x)\rt|^2 \overset{(\ref{eqs2})}=\\
			&\qd\qd\qd\frac{1}{2}\lt(\lt( \sin(\psi(x+h))-\sin(\psi(x))\rt)^2+\lt(\cos(\psi(x+h))-\cos(\psi(x))\rt)^2\rt)\\
			&\qd\qd\qd+\frac{1}{18}\lt(\lt( \sin(3\psi(x+h))-\sin(3\psi(x))\rt)^2+\lt(\cos(3\psi(x+h))-\cos(3\psi(x))\rt)^2\rt)\\
			& \qd\qd\qd =\frac{1}{2}\lt(2-2\sin(\psi(x+h))\sin(\psi(x))-2\cos(\psi(x+h))\cos(\psi(x))\rt)\\
			& \qd\qd\qd\qd +\frac{1}{18}\lt(2-2\sin(3\psi(x+h))\sin(3\psi(x))-2\cos(3\psi(x+h))\cos(3\psi(x)) \rt)\\
			&\qd\qd\qd\overset{(\ref{eqs6}),(\ref{eqs7})}{=}\lt(1-\cos(\alpha_h(x))\rt)+\frac{1}{9}\lt(1-\cos(3\alpha_h(x))\rt) \\
			&\qd\qd\qd\overset{(\ref{eqz20.5})}{=} \frac{10}{9}-\frac{2}{3}\cos(\alpha_h(x))-\frac{4}{9}\cos^3(\alpha_h(x)).
		\end{split}
	\end{equation*}}
	When $\alpha_h$ is sufficiently small, we have
	\begin{equation*}
		\frac{10}{9}-\frac{2}{3}\cos(\alpha_h(x))-\frac{4}{9}\cos^3(\alpha_h(x)) = \frac{10}{9}-\frac{2}{3}\lt(1-\frac{\alpha_h^2}{2}\rt)-\frac{4}{9}\lt(1-\frac{\alpha_h^2}{2}\rt)^3 +o(\alpha_h^2) = \alpha_h^2+o(\alpha_h^2).
	\end{equation*}
	
	\blue{\em Step 3. \rm We have
		\begin{equation*}
		\det\lt(DF(x+h)-DF(x)\rt)> c_0\lt|DF(x+h)-DF(x)\rt|^4
		\end{equation*}
		for some constant $c_0$ independent of $x$ and $h$.}
	
	\blue{\em Proof of Step 3. } \rm It follows from (\ref{eqs100}) and (\ref{eqs10}) that there exist $\delta>0$ and $c_1>0$, such that, for all $0\leq \alpha_h<\delta$, we have
	\begin{equation}\label{eqs11}
	\det\lt(DF(x+h)-DF(x)\rt)> c_1\lt|DF(x+h)-DF(x)\rt|^4.
	\end{equation}
	Let $f(t)=\frac{4}{9}-\frac{2}{3}\cos(t)+\frac{2}{9}\cos^3(t)$. Then $f'(t)=\frac{2}{3}\sin^3(t)$. Note that $f(0)=0$ and $\int_{0}^{t}\frac{2}{3}\sin^3(s)ds>0$ for all $t\in(0,2\pi)$, since $\sin^3(t)>0$ and \red{$\sin^3(t)=-\sin^3(t+\pi)$ for $t\in(0,\pi)$}. Therefore, for all $t\in(0,2\pi)$, we have $f(t)=f(0)+\int_{0}^{t}\frac{2}{3}\sin^3(s)ds >0$. By periodicity of the function $f$,  it is clear that $f(t)\geq 0$ for all $t\in \mathbb R$, and $f(t)=0$ if and only if $t=2k\pi$, $k\in\mathbb{Z}$. Similarly, given $0<\delta<\pi$, by the odd symmetry of $f'(t)$ with respect to $t=\pi$, we have $\int_{\delta}^{t}f'(s)ds\geq 0$ for all $\delta\leq t\leq 2\pi-\delta$. As a consequence, for all $\delta\leq t\leq 2\pi-\delta$, we have $f(t)\geq f(\delta)$. Thus \xblue{for all $\delta\leq \alpha_h\leq 2\pi-\delta$} we have 
		\begin{equation}\label{eqz3}
			\frac{4}{9}-\frac{2}{3}\cos(\alpha_h(x))+\frac{2}{9}\cos^3(\alpha_h(x)) \geq \frac{4}{9}-\frac{2}{3}\cos(\delta)+\frac{2}{9}\cos^3(\delta)>0.
		\end{equation}
		\xblue{Note that $\lt|DF(x+h)-DF(x)\rt|^4$ is uniformly bounded for all $x$ and $h$ such that $DF(x),DF(x+h)\in K$.} Therefore, it follows from (\ref{eqs100}) and (\ref{eqz3}) that there exists some $c_2>0$ such that \xblue{for all $\delta\leq \alpha_h\leq 2\pi-\delta$} 
		\begin{equation}
		\label{eqs12}
		\det\lt(DF(x+h)-DF(x)\rt)> c_2\lt|DF(x+h)-DF(x)\rt|^4.
		\end{equation}
	Since $f(t)$ is even with respect to $t=0$ and periodic with period $2\pi$, and so is the function $g(t)=\frac{10}{9}-\frac{2}{3}\cos(t)-\frac{4}{9}\cos^3(t)$, it is clear that the estimate (\ref{eqs11}) also holds for $2\pi-\delta<\alpha_h< 2\pi$. Combining (\ref{eqs11}) with (\ref{eqs12}), and using the periodicity of the functions $f$ and $g$, we conclude that
	\begin{equation*}
		\det\lt(DF(x+h)-DF(x)\rt)> c_0\lt|DF(x+h)-DF(x)\rt|^4,
	\end{equation*}
	where $c_0=\min\{c_1,c_2\}>0$ is independent of $x$ and $h$.
	\end{proof}

	%
	%%%%%%%%%%%%%%%%%%%%%%%%%%%%%%%%%%%%%%%%%%%%%%%%%%%%%%%%%%%%%%%%%%%%%%
	%
	% Fractional Sobolev regularity 
	%
	%%%%%%%%%%%%%%%%%%%%%%%%%%%%%%%%%%%%%%%%%%%%%%%%%%%%%%%%%%%%%%%%%%%%%%%%%
	%
	%
	%
	\subsubsection*{Proof of Theorem \ref{T2} completed.} 
	\xblue{Here we follow the idea in the proof of Theorem 3 in \cite{sverak} to show the regularity of $Dv$.} Let $\Omega'\subset \subset \Omega$ and $\gamma:=\mathrm{dist}(\Omega',\partial \Omega)>0$. Let $\eta\in C_c^{\infty}(\Omega)$ be such that $\eta\equiv 1$ on $\Omega'$ \xblue{and $\mathrm{dist}(\spt(\eta),\partial \Omega)\geq\frac{\gamma}{2}$}. \xblue{Given $e\in \mathbb{S}^1$ and $h\in\Rb$ satisfying $0<h<\frac{\gamma}{2}$, we have $DF(x),DF(x+he)\in K$ for a.e. $x\in\spt(\eta)$. It follows from Lemma \ref{l7} that}
		\begin{eqnarray}
		\label{eqzzz35}
		\det\lt(\eta(x)\frac{DF(x+he)-DF(x)}{h}\rt)&=&\frac{\eta(x)^2}{h^2}\det\lt(DF(x+he)-DF(x)\rt)\nn\\
		&\overset{(\ref{eqzz41})}{\geq}& c_0 \eta(x)^2 \frac{\lt|DF(x+he)-DF(x)\rt|^4}{h^2} \xblue{\text{ for a.e. }x\in\spt(\eta)}. 
		\end{eqnarray}
	Using the identity $\det(A+B)=\blue{\det(A)+\det(B)+A:\mathrm{Cof}(B)}$, \blue{where $\mathrm{Cof}\lt(\begin{array}{cc}   a_{11}    & a_{12} \\ a_{21}  &    a_{22}    \end{array}\rt)=\lt(\begin{array}{cc}   a_{22}    & -a_{21} \\ -a_{12}  &    a_{11}    \end{array}\rt)$, we have}
	\begin{eqnarray}
	0&=&\int_{\Omega} \det\lt(D\lt(\eta(x)\lt(\frac{F(x+he)-F(x)}{h} \rt) \rt) \rt) dx\nn\\
	&=&\int_{\Omega} \det\lt(D\eta(x)\otimes \lt(\frac{F(x+he)-F(x)}{h} \rt)+\eta(x)\lt(\frac{DF(x+he)-DF(x)}{h}  \rt) \rt) dx\nn\\
	&=&\int_{\Omega} \det\lt(D\eta(x)\otimes \lt(\frac{F(x+he)-F(x)}{h} \rt)\rt)+D\eta(x)\otimes \lt(\frac{F(x+he)-F(x)}{h} \rt):\mathrm{Cof}\lt( \eta(x)\lt(\frac{DF(x+he)-DF(x)}{h}  \rt) \rt)\nn\\
	&~&\qd\qd\qd+\det\lt( \eta(x)\lt(\frac{DF(x+he)-DF(x)}{h}  \rt)  \rt) dx.\nn
	\end{eqnarray}
	Since $\det(a\otimes b)=0$ for any $a,b\in \R^2$, the above simplifies to 
	\begin{equation}\label{eqx200}
	\begin{split}
	0&=\int_{\Omega} D\eta(x)\otimes \lt(\frac{F(x+he)-F(x)}{h} \rt):\mathrm{Cof}\lt( \eta(x)\lt(\frac{DF(x+he)-DF(x)}{h}  \rt) \rt)\\
	&\qd\qd\qd+\det\lt( \eta(x)\lt(\frac{DF(x+he)-DF(x)}{h}  \rt)  \rt) dx.
	\end{split}
	\end{equation}
	Using (\ref{eqzzz35}), (\ref{eqx200}) \xblue{and H\"{o}lder's inequality}, we have
	\begin{equation}
	\label{eqss2}
	\begin{split}
	&\int_{\Omega} \eta(x)^2 \frac{\lt|DF(x+he)-DF(x)\rt|^4}{h^2} dx\\
	&\qd\qd \overset{(\ref{eqzzz35})}{\leq} \frac{1}{c_0}\int_{\Omega} \det\lt( 
	\eta(x)\lt(\frac{DF(x+he)-DF(x)}{h}  \rt)\rt)dx\\
	&\qd\qd \overset{(\ref{eqx200})}{\leq} \frac{1}{c_0}\int_{\Omega} \frac{\lt|D\eta(x)\rt|}{\sqrt{h}} \lt|\frac{F(x+he)-F(x)}{h} \rt|\lt|\eta(x)\rt|\lt|\frac{DF(x+he)-DF(x)}{\sqrt{h}} \rt|dx\\
	&\qd\qd  \leq\frac{\blue{C(\Omega)}}{c_0} \|D\eta\|_{L^{\infty}(\Omega)}\|\sqrt{\eta}\|_{L^{\infty}(\Omega)}\mathrm{Lip}(F)\frac{1}{\sqrt{h}} \lt(\int_{\Omega} \eta(x)^2  \frac{\lt|DF(x+he)-DF(x)\rt|^4}{h^2} dx   \rt)^{\frac{1}{4}}\\
	&\qd\qd  \leq \frac{\blue{C(\Omega,\gamma)}}{\sqrt{h}}\lt(\int_{\Omega} \eta(x)^2\frac{\lt|DF(x+he)-DF(x)\rt|^4}{h^{2}} dx   \rt)^{\frac{1}{4}},
	\end{split}
	\end{equation}
	\blue{where the constant $C(\Omega,\gamma)$ depends only on $\Omega$ and $\gamma$.}

	Given $\beta\in (0,\frac{4}{3})$, it follows from (\ref{eqss2}) that
	\begin{equation}
	\label{eqzzz107}
	\begin{split}
	&\frac{1}{h^{\beta}}\int_{\Omega} \eta(x)^2\frac{\lt|DF(x+he)-DF(x)\rt|^4}{h^2} dx\\
	&\qd\qd\qd\overset{(\ref{eqss2})}{\leq}\frac{1}{h^{\beta}}\frac{\blue{C(\Omega,\gamma)}}{\sqrt{h}}\lt(\int_{\Omega} \eta(x)^2\frac{\lt|DF(x+he)-DF(x)\rt|^4}{h^{2}} dx   \rt)^{\frac{1}{4}}\\
	&\qd\qd\qd=\frac{\blue{C(\Omega,\gamma)}}{h^{\frac{1}{2}}}\frac{1}{h^{\frac{3\beta}{4}}}\lt(\int_{\Omega} \eta(x)^2\frac{\lt|DF(x+he)-DF(x)\rt|^4}{h^{2+\beta}} dx   \rt)^{\frac{1}{4}}.
	\end{split}
	\end{equation}
	\xblue{Note that the above estimate (\ref{eqzzz107}) holds for all $e\in \mathbb{S}^1$ and for all $0<h<\frac{\gamma}{2}$.} So for $0<R<\xblue{\frac{\gamma}{2}}$ we have
	\begin{eqnarray}
	\label{eqx280}
	&~&\int_{B_R} \int_{\Omega} \eta(x)^2\frac{\lt|DF(x+y)-DF(x)\rt|^4}{\lt|y\rt|^{2+\beta}} dx dy\nn\\
	&~&\qd\qd \overset{(\ref{eqzzz107})}{\leq} \blue{C(\Omega,\gamma)}\int _{B_R}\frac{1}{\lt|y\rt|^{\frac{3\beta+2}{4}}}\lt(\int_{\Omega} \eta(x)^2\frac{\lt|DF(x+y)-DF(x)\rt|^4}{\lt|y\rt|^{2+\beta}} dx   \rt)^{\frac{1}{4}} dy.
	\end{eqnarray}
	Now by Holder's inequality  
	\begin{eqnarray}
	\label{eqx29}
	&~&\int _{B_R}\frac{1}{\lt|y\rt|^{\frac{3\beta+2}{4}}}\lt(\int_{\Omega} \eta(x)^2\frac{\lt|DF(x+y)-DF(x)\rt|^4}{\lt|y\rt|^{2+\beta}} dx   \rt)^{\frac{1}{4}} dy\nn\\
	&~&\qd\qd\leq \lt(\int_{B_R} \frac{1}{\lt|y\rt|^{\beta+\frac{2}{3}}}  dy\rt)^{\frac{3}{4}}
	\lt(\int_{B_R} \int_{\Omega} \eta(x)^2\frac{\lt|DF(x+y)-DF(x)\rt|^4}{\lt|y\rt|^{2+\beta}} dx dy    \rt)^{\frac{1}{4}}.
	\end{eqnarray}
	As $\beta\in (0,\frac{4}{3})$, let $\delta:=2-(\beta+\frac{2}{3})\blue{>0}$, then $\beta+\frac{2}{3}=2-\delta$. We have
	\begin{equation}
	\label{eqx210}
	\int_{B_R} \frac{1}{\lt|y\rt|^{2-\delta}}  dy=2\pi\int_{0}^R\frac{1}{r^{2-\delta}} r dr=2\pi \int_{0}^R r^{-1+\delta} dr=\frac{2\pi}{\delta} R^{\delta}. 
	\end{equation}
	Putting this together with (\ref{eqx280})-(\ref{eqx29}) we have 
	\begin{eqnarray}
	\label{eqx28}
	&~&\int_{B_R} \int_{\Omega} \eta(x)^2\frac{\lt|DF(x+y)-DF(x)\rt|^4}{\lt|y\rt|^{2+\beta}} dx dy\nn\\
	&~&\qd\qd \overset{(\ref{eqx280}),(\ref{eqx29}),(\ref{eqx210})}{\leq} \blue{C(\Omega,\gamma)}\lt( \frac{2\pi}{\delta} R^{\delta} \rt)^{\frac{3}{4}} \lt(\int_{B_R} \int_{\Omega} \eta(x)^2\frac{\lt|DF(x+y)-DF(x)\rt|^4}{\lt|y\rt|^{2+\beta}} dx dy    \rt)^{\frac{1}{4}}.
	\end{eqnarray}
	Thus, \xblue{noting $\delta=2-(\beta+\frac{2}{3})>0$}, we deduce from (\ref{eqx28}) that
	\begin{equation}
	\label{eqx28.6}
	\int_{B_R} \int_{\Omega} \eta(x)^2\frac{\lt|DF(x+y)-DF(x)\rt|^4}{\lt|y\rt|^{2+\beta}} dx dy<\blue{C(\Omega,\gamma,\beta)}
	\end{equation}
	for some constant $C(\Omega,\gamma,\beta)$ depending only on $\Omega$, $\gamma$ and $\beta$. Note that $\eta(x)\equiv 1$ for $x\in \Omega'$. Therefore, we deduce from (\ref{eqx28.6}) that
		\begin{equation*}
		\int_{B_R} \int_{\Omega'} \frac{\lt|DF(x+y)-DF(x)\rt|^4}{\lt|y\rt|^{2+\beta}} dx dy<\blue{C(\Omega,\gamma,\beta)}.
		\end{equation*}
		It follows that
		\begin{equation*}
		\begin{split}
		&\int_{\Omega'} \int_{\Omega'} \frac{\lt|DF(x)-DF(w)\rt|^4}{\lt|x-w\rt|^{2+\beta}} dw dx\\
		&\qd\qd \leq \int_{\Omega'}\int_{B_R(x)} \frac{\lt|DF(x)-DF(w)\rt|^4}{\lt|x-w\rt|^{2+\beta}} dw dx + \int_{\Omega'}\int_{\Omega'\setminus B_R(x)} \frac{\lt|DF(x)-DF(w)\rt|^4}{\lt|x-w\rt|^{2+\beta}} dw dx\\
		&\qd\qd < C(\Omega,\gamma,\beta) + \frac{1}{R^{2+\beta}}\int_{\Omega'}\int_{\Omega'\setminus B_R(x)} \lt|DF(x)-DF(w)\rt|^4 dw dx < C.
		\end{split}
		\end{equation*}
	So this implies $DF\in \blue{W^{\frac{\beta}{4},4}(\Omega')}$. \xblue{Recall that $F(x_1,x_2)=\lt(\mathrm{Re}\lt(v(x_1+ix_2)\rt),\mathrm{Im}\lt(v(x_1+ix_2)\rt)\rt)$}. Therefore we have established (\ref{eqy60.4}).

	%%%%%%%%%%%%%%%%%%%%%%%%%%%%%%%%%%%%%
	
	Now from (\ref{eqss2}) we have that for any $y\in B_{\xblue{\frac{\gamma}{2}}}(0)$
	\begin{equation*}
	\int_{\Omega'} \frac{\lt|DF(x+y)-DF(x)\rt|^4}{|y|^2} dx \leq \int_{\Omega} \eta(x)^2 \frac{\lt|DF(x+y)-DF(x)\rt|^4}{|y|^2} dx \overset{(\ref{eqss2})}{\leq} \lt(\frac{\blue{C(\Omega,\gamma)}}{\sqrt{|y|}}\rt)^{\frac{4}{3}}.
	\end{equation*}
	It follows that
	\begin{equation}\label{eq201}
	\int_{\Omega'} \frac{\lt|DF(x+y)-DF(x)\rt|^4}{|y|^{\frac{4}{3}}} dx \leq \blue{\wt C(\Omega,\gamma)}.
	\end{equation}
	\xblue{Given $0<\ep<\frac{\gamma}{2}$,} integrating the above with respect to $y$ over $B_{\ep}(0)$ and using the fact that $|y|\leq \ep$ for all $y\in B_{\ep}(0)$, we obtain
	\begin{equation*}
	\begin{split}
	&\int_{\Omega'}\int_{B_{\ep}(0)} \frac{\lt|DF(x+y)-DF(x)\rt|^4}{\ep^{2+\frac{4}{3}}} dydx\\
	&\qd\qd \leq \frac{1}{\ep^2}\int_{\Omega'}\int_{B_{\ep}(0)} \frac{\lt|DF(x+y)-DF(x)\rt|^4}{|y|^{\frac{4}{3}}} dydx\\
	&\qd\qd \overset{(\ref{eq201})}{\leq} \frac{\blue{\wt C(\Omega,\gamma)}}{\ep^2} \int_{B_{\ep}(0)}1dy = \pi\blue{\wt C(\Omega,\gamma)}.
	\end{split}
	\end{equation*}
	This establishes (\ref{eqss1}). $\Box$ \nl
	
	As a corollary of Theorem \ref{T2}, we have
	
	\begin{a6}\label{c5}
		Let $\Omega \subset \R^2$ be a \bblue{bounded simply-connected} domain and $u \in E(\Omega)$, \bblue{where $E(\Omega)$ is defined in (\ref{eqz4})}. Then \blue{$\nabla u \in W^{\xblue{\sigma},4}_{loc}(\Omega)$ for all $0<\xblue{\sigma} < \frac{1}{3}$}. Further, for any $\Omega'\subset\subset \Omega$, there exists a constant $C$ such that 
		\begin{equation}\label{c5.0}
		\int_{\Omega'}\int_{B_{\ep}\bblue{(0)}}\frac{\lt|\nabla u(x+y)-\nabla u(x)\rt|^4}{\ep^{2+\frac{4}{3}}} dy dx<C
		\end{equation}
		for all $\ep$ sufficiently small, \blue{where the above constant $C$ is independent of $\ep$}.
	\end{a6}
	
	\begin{proof}[Proof of Corollary \ref{c5}] 
	Since $u\in E(\Omega)$, it follows from Theorem \ref{T3} that there exists $F_u$ such that 
	\begin{equation*}
	DF_u=\lt(\begin{array}{cc} u_{,2}\lt(1-u_{,1}^2-\frac{u_{,2}^2}{3}\rt)  & u_{,1}\lt(1-u_{,2}^2-\frac{u_{,1}^2}{3}\rt)  \\ -u_{,1}\lt(1-\frac{2 u_{,1}^2}{3}\rt)  &  u_{,2}\lt(1-\frac{2 u_{,2}^2}{3}\rt) \end{array}\rt)
	\end{equation*}
	and therefore $DF_u \in K$ a.e. in $\Omega$, where the space $K$ is defined in (\ref{eqy29}). Using (\ref{eqz59}) we have that 
	\begin{equation}\label{c5.2}
	u_{,1}= F_{u,2}^1 - F_{u,1}^2 \qd\text{ and }\qd u_{,2}= F_{u,1}^1 + F_{u,2}^2 \qd\text{ a.e. in } \Omega.
	\end{equation}
	From Theorem \ref{T2}, we have \blue{$DF_u \in W^{\xblue{\sigma},4}_{loc}(\Omega)$ for all $\xblue{\sigma} < \frac{1}{3}$}, and for any $\Omega'\subset\subset\Omega$, 
	\begin{equation}\label{c5.3}
	\int_{\Omega'}\int_{B_{\ep}\bblue{(0)}}\frac{\lt|DF_u(x+y)-DF_u(x)\rt|^4}{\ep^{2+\frac{4}{3}}} dy dx<C
	\end{equation}
	for some constant $C$ independent of $\ep$. It follows from \eqref{c5.2} that \blue{$\nabla u \in W^{\xblue{\sigma},4}_{loc}(\Omega)$ for all $\xblue{\sigma} < \frac{1}{3}$}. 
	
	By the inequality $\lt|A+B\rt|^4\leq 8(\lt|A\rt|^4+\lt|B\rt|^4)$ we have that 
	\begin{eqnarray}
	\label{eqzzz108}
	\lt|u_{,1}(x+y)-u_{,1}(x)\rt|^4&\overset{(\ref{c5.2})}{=}&\lt|\lt(F^1_{u,2}(x+y)-F^2_{u,1}(x+y)\rt)-\lt(F^1_{u,2}(x)-F^2_{u,1}(x)\rt)\rt|^4\nn\\
	&\leq& 8\lt|F^1_{u,2}(x+y)-F^1_{u,2}(x)\rt|^4+8\lt| F^2_{u,1}(x+y)-F^2_{u,1}(x)\rt|^4\nn\\
	&\leq& C\lt|DF_{u}(x+y)- DF_{u}(x)\rt|^4.
	\end{eqnarray}
	In the same way we can show that
	\begin{equation}
	\label{eqzzz109}
	\blue{\lt|u_{,2}(x+y)-u_{,2}(x)\rt|\leq C\lt|DF_{u}(x+y)- DF_{u}(x)\rt|^4}.
	\end{equation}
	Thus 
	\begin{equation}\label{c5.4}
	\lt|\nabla u(x+y)-\nabla u(x)\rt|^4 \overset{(\ref{eqzzz108}),(\ref{eqzzz109})}{\leq} C\lt|DF_u(x+y)-DF_u(x)\rt|^4
	\end{equation}
	for some pure constant $C$. Finally, putting \eqref{c5.3} and \eqref{c5.4} together, we immediately obtain \eqref{c5.0}. 
	\end{proof}

	%%%%%%%%%%%%%%%%%%%%%%%%%%%%%%%%%%%%%%%%%%%%%%%%%%%%%%%
	%
	%   Vanishing of a class of entropies
	%
	%%%%%%%%%%%%%%%%%%%%%%%%%%%%%%%%%%%%%%%%%%%%%%%%%%%%%%%%%%%%

	\section{Vanishing of \bblue{harmonic} entropies}
	\bblue{Recall the definition of entropies in (\ref{eqi151})}. We first recall a few lemmas from \cite{mul2}.
	
	\begin{a1}[\cite{mul2}, Lemma 1]\label{l12}
		Let $\Phi\in C_{c}^{\infty}(\Rb^2;\R^2)$ be an entropy. Then there exists a $\Psi\in C_{c}^{\infty}(\Rb^2;\Rb^2)$ such that
		\begin{equation}\label{l12.2}
		D\Phi(z)+2\Psi(z)\otimes z \quad \text{ is isotropic for all } z.
		\end{equation}
		Consequently, we have
		\begin{equation}\label{l12.1}
		\Psi_1(z)=-\frac{1}{2z_2}\Phi_{1,2}(z) \quad \text{ and } \quad \Psi_{2}(z)=-\frac{1}{2z_1}\Phi_{2,1}(z).
		\end{equation}
	\end{a1}
	
	\begin{a1}[\cite{mul2}, Lemma 2]\label{l18}
		Let $\Phi\in C_{c}^{\infty}(\Rb^2;\Rb^2)$ and $\Psi \in C_{c}^{\infty}(\Rb^2;\Rb^2)$ satisfy \eqref{l12.2}. Let $m \in H^1(\Omega;\Rb^2)$ satisfy
		\begin{equation*}
			\na \cdot m = 0 \quad \text{a.e. in } \Omega.
		\end{equation*}
		Then
		\begin{equation}\label{l18.1}
		\na \cdot \lt[\Phi(m)\rt] = \Psi(m)\cdot \na\lt(1-\lt|m\rt|^2\rt) \quad \text{a.e. in } \Omega.
		\end{equation}
	\end{a1}
	
	\begin{a1}[\cite{mul2}, Lemma 3]\label{l13}
		There is a one-to-one correspondence between entropies $\Phi\in C_{c}^{\infty}(\Rb^2;\Rb^2)$ and functions $\varphi\in C_{c}^{\infty}(\Rb^2)$ with $\varphi(0)=0$ via
		\begin{equation}\label{s4.1}
		\Phi(z)=\varphi(z)z+\lt(\nabla\varphi(z)\cdot z^{\perp}\rt)z^{\perp},
		\end{equation} 
		where $z^{\perp}=(-z_2,z_1)$ is the anticlockwise rotation of $z$ by $\frac{\pi}{2}$.
	\end{a1}
	
	Using the above lemmas, we have the following relationship between $\Psi$ and $\varphi$.
	
	\begin{a1}\label{l17}
		Let $\Phi\in C_{c}^{\infty}(\Rb^2;\Rb^2)$ be an entropy, and $\Psi$ and $\varphi$ be the functions related to $\Phi$ through Lemmas \ref{l12} and \ref{l13}, respectively. Then we have
		\begin{equation}\label{s4.2}
		\Psi_{1,2}(z)-\Psi_{2,1}(z)=\frac{1}{2}\nabla\lt(\Delta\varphi\rt)\cdot z^{\perp}.
		\end{equation}
	\end{a1}
	
	\begin{proof}
		Using formula \eqref{s4.1}, we have
		\begin{equation*}
			\Phi_{1,2}(z)=2z_2\varphi_{,1}(z)+z_2^2\varphi_{,12}(z)-z_1z_2\varphi_{,22}(z),\quad \Phi_{2,1}(z)=2z_1\varphi_{,2}(z)+z_1^2\varphi_{,12}(z)-z_1z_2\varphi_{,11}(z).
		\end{equation*}
		Putting the above into \eqref{l12.1}, we obtain
		\begin{equation*}
			\Psi_1(z)=-\varphi_{,1}(z)-\frac{z_2}{2}\varphi_{,12}(z)+\frac{z_1}{2}\varphi_{,22}(z),\quad \Psi_2(z)=-\varphi_{,2}(z)-\frac{z_1}{2}\varphi_{,12}(z)+\frac{z_2}{2}\varphi_{,11}(z).
		\end{equation*}
		By direct calculations, we have
		\begin{equation*}
			\Psi_{1,2}(z)=-\frac{3}{2}\varphi_{,12}(z)-\frac{z_2}{2}\varphi_{,122}(z)+\frac{z_1}{2}\varphi_{,222}(z),\quad \Psi_{2,1}(z)=-\frac{3}{2}\varphi_{,12}(z)-\frac{z_1}{2}\varphi_{,112}(z)+\frac{z_2}{2}\varphi_{,111}(z).
		\end{equation*}
		Hence, we have
		\begin{equation*}
			\Psi_{1,2}(z)-\Psi_{2,1}(z)=\frac{1}{2}\lt(\varphi_{,122}+\varphi_{,111},\varphi_{,112}+\varphi_{,222}\rt)\cdot(-z_2,z_1)=\frac{1}{2}\nabla\lt(\Delta\varphi\rt)\cdot z^{\perp}.
		\end{equation*}
	\end{proof}
	
	Given a function $u \in W^{1,\bblue{\infty}}(\Omega)$, for all $\ep>0$, we denote $u_{\ep} = u\ast\rho_{\ep}$, where $\rho_{\ep}$ is the standard \bblue{convolution kernel} supported in $B_{\ep}(0)\subset \Rb^2$. Very often in this paper, we use the following notation
	\begin{equation}\label{s4.3}
	w:= \lt(\nabla u\rt)^{\perp} = \lt(-u_{,2},u_{,1}\rt).
	\end{equation}
	Then we have
	\begin{equation}\label{s4.4}
	w_{\ep}=\lt( w_{\ep}^1,w_{\ep}^2 \rt) = \lt( -u_{\ep,2}, u_{\ep,1} \rt).
	\end{equation}
	\green{The main result of this section is the following theorem.}
	
	\begin{a2}\label{T4}
		Let $\Omega \subset \R^2$ be a bounded \bblue{simply-connected} domain. Let $\Phi\in C_{c}^{\infty}(\Rb^2;\Rb^2)$ be an entropy, and $\varphi\in C_{c}^{\infty}(\Rb^2)$ with $\varphi(0)=0$ be the smooth function related to $\Phi$ through \eqref{s4.1}. In addition, we assume that 
		\begin{equation}\label{T4.4}
		\nabla\lt(\Delta\varphi\rt)\cdot z^{\perp}=0 \qd\text{ for all } z\in\R^2.
		\end{equation}
		Then, for all $u \in E(\Omega)$, \bblue{where $E(\Omega)$ is defined in (\ref{eqz4}),} we have 
		\begin{equation*}
		\nabla\cdot\lt[\Phi(\na u^{\perp})\rt]=0
		\end{equation*}
		in the sense of distributions.
	\end{a2}
	
	\begin{proof}[Proof of Theorem \ref{T4}]
		
		\xblue{Given $\Omega'\subset\subset\Omega$,} let $\zeta\in C_c^{\infty}(\Omega')$ be a test function. Recall $w\overset{(\ref{s4.3})}{=}(\na u)^{\perp}=(-u_{,2},u_{,1})$. So as in Step 6 of the proof of Proposition 3 \cite{DeI} we have 
		\begin{eqnarray}\label{eqs72}
		&~&\int_{\Omega'} \zeta(x) \nabla \cdot \lt[ \Phi(w_{\ep})\rt] dx\nn\\
		&~&\qd\qd\overset{(\ref{l18.1})}{=}\int_{\Omega'} \zeta(x) \Psi(w_{\ep})\cdot \na \lt(1-\lt|w_{\ep}\rt|^2 \rt) dx\nn\\
		&~&\qd\qd=\overbrace{\int_{\Omega'} \zeta(x) \na \cdot \lt[\Psi(w_{\ep}) \lt(1-\lt|w_{\ep}\rt|^2 \rt)  \rt] dx}^{I_{\ep}}
		-\overbrace{\int_{\Omega'} \zeta(x)  \lt(1-\lt|w_{\ep}\rt|^2 \rt) \na\cdot \lt[ \Psi(w_{\ep}) \rt] dx}^{\Pi_{\ep}}.
		\end{eqnarray}
		Since $\Psi(w_{\ep}) \lt(1-\lt|w_{\ep}\rt|^2 \rt)\overset{L^1}{\rightarrow} 0$, integrating by parts we see that 
		\begin{equation}
		\label{eqs511}
		I_{\ep}\rightarrow 0. 
		\end{equation}
		In the following, we show that, under the additional assumption \eqref{T4.4}, we have
		\begin{equation*}
		\Pi_{\ep}\rightarrow 0.
		\end{equation*}
		Thus, \bblue{we have}
$$
\xred{\int_{\Omega'} \Phi(w)\cdot \na \zeta dx=\lim_{\ep\rightarrow 0} \int_{\Omega'} \Phi(w_{\ep})\cdot \na \zeta dx= -\lim_{\ep\rightarrow 0} \int_{\Omega'} \na \cdot \lt[\Phi(w_{\ep})\rt] \zeta dx=0},
$$		
\bblue{from which Theorem \ref{T4} will follow}.

		% % % % % % % % % % % % % % % % % % % % % % % % % % % % % % % % % % % % %
		
		We will need several lemmas. First, we provide the following lemma, which will be used repeatedly.
		
		\begin{a1}\label{l11}
			\bblue{Let $\Omega$ be as in Theorem \ref{T4} and $u\in E(\Omega)$.} Given $\Omega'\subset\subset \Omega$, there exists a constant $\ep_{0}=\ep_{0}(\Omega')$ such that, for all $\ep<\ep_{0}$, and for all $r\geq 4$ and $f \in L^{r}(\Omega')$, we have
			\begin{equation}\label{l114}
			\int_{\Omega'}\lt|1-|\nabla u_{\ep}|^2\rt|\lt|u_{\ep,mn}\rt| \lt|f\rt| dx \leq C \lVert f \rVert_{L^{r}(\Omega')}
			\end{equation}
			for all $m=1, 2$ and $n=1, 2$, and for some constant $C$ independent of $\ep$. 
			
			Consequently, if $g_j\rightarrow g$ in $L^r(\Omega)$, then for any sequence $\{\ep_{j}\}$ such that $0<\ep_{j}<\ep_{0}$ for all $j$, we have 
			\begin{equation}\label{l111}
			\int_{\Omega'}\lt|1-|\nabla u_{\ep_j}|^2\rt|\lt|u_{\ep_j,mn}\rt|\lt|g-g_j\rt| dx\rightarrow 0 \quad \text{ as } j\rightarrow \infty
			\end{equation}
			for all $m=1, 2$ and $n=1, 2$.
		\end{a1}
		
		\begin{proof}
			Given $r\geq 4$ and $f \in L^{r}(\Omega')$, by H\"{o}lder's inequality, we have
			\begin{equation}\label{l112}
			\int_{\Omega'}\lt|1-|\nabla u_{\ep}|^2\rt| \lt|u_{\ep,mn}\rt| \lt|f\rt| dx\leq \lt(\int_{\Omega'}\lt|1-|\nabla u_{\ep}|^2\rt|^{r'}\lt|u_{\ep,mn}\rt|^{r'}dx\rt)^{\frac{1}{r'}} \lVert f \rVert_{L^r(\Omega')},
			\end{equation}
			where $\frac{1}{r}+\frac{1}{r'}=1$. Now as in (i) and (ii) of Step 6 of the proof of Proposition 3 \cite{DeI}, we have 
			that
			\begin{equation}
			\label{eqx305}
			1-\lt|w_{\ep}(x)\rt|^2\leq \frac{2 \|\rho\|_{L^{\infty}(\R^2)}}{\ep^2}\int_{B_{\ep}} \lt|w(x-z)-w(x)\rt|^2 dz,
			\end{equation} 
			and
			\begin{equation}
			\label{eqx306}
			\lt|w_{\ep,j}(x)\rt|\leq \frac{\|\na \rho\|_{L^{\infty}(\R^2)}}{\ep^3}\int_{B_{\ep}} \lt|w(x-z)-w(x)\rt| dz,
			\end{equation}
			where recall that we defined the vector fields $w$ and $w_{\ep}$ in \eqref{s4.3} and \eqref{s4.4}, respectively. For the convenience of the reader, we take the proofs of \eqref{eqx305}-\eqref{eqx306} in \cite{DeI} and put them into Lemma \ref{l21} in \xblue{the Appendix}.
			
			Note that since $r\geq 4$, we have $r'\leq \frac{4}{3}$ and, therefore, $\frac{3r'}{4}\leq 1$. Now arguing very similarly to (iii) of Step 6 of  Proposition 3 \cite{DeI}, we have
			\begin{equation}\label{l113}
			\begin{split}
			&\int_{\Omega'} \lt|1-\lt|\nabla u_{\ep}\rt|^2\rt|^{r'}\lt|u_{\ep,mn}\rt|^{r'}dx\\
			\overset{\xred{(\ref{eqx305}),(\ref{eqx306})}}{\leq} & \frac{C}{\epsilon^{r'}}\int_{\Omega'}\lt(\Xint{-}_{B_{\ep}}\lt|w(x-z)-w(x)\rt|^2dz\rt)^{r'}\lt(\Xint{-}_{B_{\ep}}\lt|w(x-z)-w(x)\rt|dz\rt)^{r'}dx\\
			\leq & \frac{C}{\epsilon^{r'}}\int_{\Omega'}\lt(\Xint{-}_{B_{\ep}}\lt|w(x-z)-w(x)\rt|^4dz\rt)^{\frac{r'}{2}}\lt(\Xint{-}_{B_{\ep}}\lt|w(x-z)-w(x)\rt|^4dz\rt)^{\frac{r'}{4}}dx\\
			= & \frac{C}{\epsilon^{r'}}\int_{\Omega'}\lt(\Xint{-}_{B_{\ep}}\lt|w(x-z)-w(x)\rt|^4dz\rt)^{\frac{3r'}{4}}dx\\
			\leq & \frac{C}{\epsilon^{r'}}\lt(\int_{\Omega'}\Xint{-}_{B_{\ep}}\lt|w(x-z)-w(x)\rt|^4dzdx\rt)^{\frac{3r'}{4}}\\
			= & C\lt(\int_{\Omega'}\int_{B_{\ep}}\frac{\lt|w(x-z)-w(x)\rt|^4}{\ep^{2+\frac{4}{3}}}dzdx\rt)^{\frac{3r'}{4}} \overset{(\ref{c5.0})}{\leq} C.
			\end{split}
			\end{equation}
			Putting \eqref{l113} into \eqref{l112}, we immediately obtain \eqref{l114}. The estimate \eqref{l111} is a direct consequence of \eqref{l114}.
		\end{proof}
		
		\begin{a1}\label{l15}
		
		\bblue{Let $\Omega$ be as in Theorem \ref{T4} and $u\in E(\Omega)$. Denote $w:=(\na u)^{\perp}$.} Given $\Omega'\subset \subset \Omega$, for all $\zeta\in C_{c}^{\infty}(\Omega')$, we have
			\begin{equation}\label{l15.1}
			\int_{\Omega'} (1-\lt|w_{\ep}(x)\rt|^2)w^1_{\ep,1}(x)\zeta(x) dx\rightarrow 0, \quad \int_{\Omega'} (1-\lt|w_{\ep}(x)\rt|^2)w^2_{\ep,2}(x)\zeta(x) dx\rightarrow 0,
			\end{equation}
			and
			\begin{equation}\label{l15.2}
			\int_{\Omega'} (1-\lt|w_{\ep}(x)\rt|^2)\lt(w^2_{\ep,1}(x)+w^1_{\ep,2}(x)\rt)\zeta(x) dx\rightarrow 0
			\end{equation}
			as $\ep \rightarrow 0$.
		\end{a1}
		
		\begin{proof}
			Given a smooth function $v$, by direct calculations, we have
			\begin{equation}\label{l15.3}
			\nabla\cdot\lt[\Sigma_{e_1 e_2}v\rt]\overset{\xred{(\ref{eqi501})}}{=}\lt(v_{,11}-v_{,22}\rt)\lt(1-|\nabla v|^2\rt)
			\end{equation}
			and
			\begin{equation}\label{l15.4}
			\nabla\cdot\lt[\Sigma_{\ep_1 \ep_2}v\rt]\overset{\xred{(\ref{eqi502})}}{=}2v_{,12}\lt(1-|\nabla v|^2\rt).
			\end{equation}
			Recall the definition of $w_{\ep}$ in \eqref{s4.4}. In particular, we have 
			\begin{equation}\label{l15.9}
			w^{1}_{\ep,1}=-u_{\ep,12}, \quad w^{2}_{\ep,2}=u_{\ep,12}, \quad w^2_{\ep,1}+w^1_{\ep,2}=u_{\ep,11}-u_{\ep,22}.
			\end{equation}
			Thus, using \eqref{s4.4}, \eqref{l15.4} and \eqref{l15.9}, we have
			\begin{equation}
			\label{l15.5}
			\begin{split}
			\int_{\Omega'} \lt(1-\lt|w_{\ep}(x)\rt|^2\rt)w^1_{\ep,1}(x)\zeta(x) dx&\overset{(\ref{s4.4})(\ref{l15.9})}{=}-\int_{\Omega'} \lt(1-\lt|\nabla u_{\ep}(x)\rt|^2\rt)u_{\ep,12}(x)\zeta \; dx\\
			&\overset{\xred{(\ref{eqi502})},(\ref{l15.4})}{=}-\frac{1}{2}\int_{\Omega'}\nabla\cdot\lt(u_{\ep,2}\lt(1-\frac{2 u_{\ep,2}^2}{3}\rt),  u_{\ep,1}\lt(1-\frac{2 u_{\ep,1}^2}{3}\rt)\rt) \zeta \;dx\\
			&=\frac{1}{2}\int_{\Omega'}\lt(u_{\ep,2}\lt(1-\frac{2 u_{\ep,2}^2}{3}\rt), u_{\ep,1}\lt(1-\frac{2 u_{\ep,1}^2}{3}\rt)\rt) \cdot\nabla\zeta\; dx.
			\end{split}
			\end{equation}
			Since $u \in W^{1,\bblue{\infty}}(\Omega)$, it follows from \eqref{l15.5} and \eqref{eqi46} that
			\begin{equation}\label{l15.7}
			\begin{split}
			&\int_{\Omega'} \lt(1-\lt|w_{\ep}(x)\rt|^2\rt)w^1_{\ep,1}(x)\zeta(x) dx\\
			&\quad\quad\overset{(\ref{l15.5})}{\rightarrow} \frac{1}{2}\int_{\Omega'}\lt(u_{,2}\lt(1-\frac{2 u_{,2}^2}{3}\rt), u_{,1}\lt(1-\frac{2 u_{,1}^2}{3}\rt)\rt) \cdot\nabla\zeta \; dx = \frac{1}{2}\int_{\Omega'} \Sigma_{\ep_1 \ep_2}u(x)\cdot \nabla\zeta(x) dx\overset{(\ref{eqi46})}{=}0.
			\end{split}
			\end{equation}
			Similarly, as $w^{2}_{\ep,2}=u_{\ep,12}$, we have
			\begin{equation*}
			\int_{\Omega'} \lt(1-\lt|w_{\ep}(x)\rt|^2\rt)w^2_{\ep,2}(x)\zeta(x) dx \rightarrow 0.
			\end{equation*}
			Next, using \eqref{s4.4}, \eqref{l15.3} and \eqref{l15.9}, we have
			\begin{equation*} %\label{l15.6}
			\begin{split}
			\int_{\Omega'} \lt(1-\lt|w_{\ep}(x)\rt|^2\rt)&\lt(\red{w^2_{\ep,1}(x)+w^1_{\ep,2}(x)}\rt)\zeta(x) dx\\
			&\overset{(\ref{s4.4})(\ref{l15.9})}{=}\int_{\Omega'} \lt(1-\lt|\nabla u_{\ep}(x)\rt|^2\rt)\lt(u_{\ep,11}(x)-u_{\ep,22}(x)\rt)\zeta(x) dx\\
			&\overset{\xred{(\ref{eqi501})},(\ref{l15.3})}{=}\frac{1}{2}\int_{\Omega'}\nabla\cdot\lt(u_{\ep,1}\lt(1-u_{\ep,2}^2-\frac{u_{\ep,1}^2}{3}\rt),  -u_{\ep,2}\lt(1-u_{\ep,1}^2-\frac{u_{\ep,2}^2}{3}\rt)\rt) \zeta \;dx.
			\end{split}
			\end{equation*}
			By (\ref{eqi46}) and the same arguments as in \eqref{l15.7}, we conclude \eqref{l15.2}.
		\end{proof}
		
		\begin{a1}\label{l14}
			\bblue{Let $\Omega$ be as in Theorem \ref{T4} and $u\in E(\Omega)$. Denote $w:=(\na u)^{\perp}$.} \xred{Given $\Omega'\subset \subset \Omega$ and} any $\red{F} \in C^{\blue{\infty}}_{c}(\Rb^2)$, we have
			\begin{equation}\label{l14.1}
			\int_{\Omega'}\lt(1-|w_{\ep}|^2\rt)w_{\ep,m}^n\lt(\red{F}(w)-\red{F}(w_{\delta})\rt)dx\rightarrow 0 \quad \text{ as } \ep,\delta\rightarrow 0
			\end{equation}
			for all $m=1, 2$ and $n=1, 2$. As a consequence, we have
			\begin{equation}\label{l14.2}
			\int_{\Omega'}\lt(1-|w_{\ep}|^2\rt)w_{\ep,m}^n\lt(\red{F}(w_{\ep})-\red{F}(w_{\delta})\rt)dx\rightarrow 0 \quad \text{ as } \ep,\delta\rightarrow 0.
			\end{equation}
		\end{a1}
		
		\begin{proof}
			First, it is clear that \red{by applying Lemma \ref{l11}, we have}
			\begin{equation}\label{l14.3}
			\begin{split}
			&\lt|\int_{\Omega'} (1-\lt|w_{\ep}(x)\rt|^2)\lt(\red{F}(w(x))-\red{F}(w_{\delta}(x))\rt)w^n_{\ep,m}(x) dx\rt|\\
			\leq& \sup_{\mathbb R^2}|D\red{F}|\int_{\Omega'} (1-\lt|w_{\ep}(x)\rt|^2)\lt|w^n_{\ep,m}(x)\rt|\lt|w(x)-w_{\delta}(x)\rt| dx\\
			\overset{(\ref{l114})}{\leq} & \red{C \sup_{\R^2} \lt|D F\rt| \|w-w_{\delta}\|_{L^r(\Omega')}}
			\end{split}
			\end{equation}
			for all $r\geq 4$. \red{ Now as $w\in L^{\infty}(\Omega)\subset L^{r}(\Omega)$ so $w_{\delta}\overset{L^{r}(\Omega)}{\rightarrow}w$. Applying this to (\ref{l14.3}), equation (\ref{l14.1}) follows.}
			
			To show \eqref{l14.2}, we write
			\begin{equation*}
				\begin{split}
					&\int_{\Omega'}\lt(1-|w_{\ep}|^2\rt)w_{\ep,m}^n\lt(\red{F}(w_{\ep})-\red{F}(w_{\delta})\rt)dx\\
					&\quad\quad=\int_{\Omega'}\lt(1-|w_{\ep}|^2\rt)w_{\ep,m}^n\lt(\red{F}(w_{\ep})-\red{F}(w)\rt)dx+\int_{\Omega'}\lt(1-|w_{\ep}|^2\rt)w_{\ep,m}^n\lt(\red{F}(w)-\red{F}(w_{\delta})\rt)dx.
				\end{split}
			\end{equation*}
			By applying \eqref{l14.1} to the above two terms on the right side, we obtain \eqref{l14.2}.
		\end{proof}

		\begin{a1}\label{L6}
			\bblue{Let $\Omega$ be as in Theorem \ref{T4} and $u\in E(\Omega)$. Denote $w:=(\na u)^{\perp}$.} \blue{Given $\Omega'\subset \subset \Omega$, for any $F \in C^{\infty}_{c}(\Rb^2)$ and any $\zeta\in C^{\infty}_c(\Omega')$, We have
				\begin{equation}\label{L6.3}
				\int_{\Omega'} \lt(1-\lt|w_{\ep}(x)\rt|^2\rt)w^1_{\ep,1}(x)F(w_{\ep}(x)) \zeta(x) dx \rightarrow 0,
				\end{equation}
				\begin{equation}\label{L6.4}
				\int_{\Omega'} \lt(1-\lt|w_{\ep}(x)\rt|^2\rt)w^2_{\ep,2}(x)F(w_{\ep}(x)) \zeta(x) dx \rightarrow 0,
				\end{equation}
				and
				\begin{equation}\label{L7.4}
				\int_{\Omega'}\lt(1-|w_{\ep}|^2\rt)\lt(w_{\ep,1}^2+w_{\ep,2}^1\rt)F(w_{\ep}(x)) \zeta(x) dx \rightarrow 0
				\end{equation}
				as $\ep\rightarrow 0$.}
		\end{a1}
		
		\begin{proof}
			We write
			\begin{equation}\label{L6.1}
			\begin{split}
			&\int_{\Omega'} \lt(1-\lt|w_{\ep}(x)\rt|^2\rt)F(w_{\ep}(x))w^1_{\ep,1}(x) \red{\zeta(x)} dx\\
			=&\int_{\Omega'} \lt(1-\lt|w_{\ep}(x)\rt|^2\rt)F(w_{\delta}(x))w^1_{\ep,1}(x) \red{\zeta(x)} dx \nn\\
			&\qd + \int_{\Omega'} \lt(1-\lt|w_{\ep}(x)\rt|^2\rt)\lt(F(w_{\ep}(x))-F(w_{\delta}(x))\rt)w^1_{\ep,1}(x) \red{\zeta(x)} dx,
			\end{split}
			\end{equation}
			where $w_{\delta}=\rho_{\delta}*w$. It follows from Lemma \ref{l15} (\ref{l15.1}) that, for any fixed $\delta>0$,
			\begin{equation}\label{L6.2}
			\int_{\Omega'} (1-\lt|w_{\ep}(x)\rt|^2)F(w_{\delta}(x))w^1_{\ep,1}(x)\red{\zeta(x)} dx\rightarrow 0 \quad \text{ as } \ep\rightarrow 0. 
			\end{equation}
			On the other hand, we obtain from Lemma \ref{l14} that
			\begin{equation}\label{L6.5}
			\int_{\Omega'} (1-\lt|w_{\ep}(x)\rt|^2)\lt(F(w_{\ep}(x))-F(w_{\delta}(x))\rt)w^1_{\ep,1}(x)\red{\zeta(x)}  dx\rightarrow 0 \quad \text{ as } \ep,\delta \rightarrow 0.
			\end{equation}
			
			Given $\alpha>0$, it follows from \eqref{L6.5} that there exist $\delta_{0}=\delta_{0}(\alpha),\ep_{0}=\ep_{0}(\alpha)>0$ sufficiently small such that, for all $\ep<\ep_{0}$, we have
			\begin{equation}\label{L6.6}
			\lt|\int_{\Omega'} (1-\lt|w_{\ep}(x)\rt|^2)\lt(F(w_{\ep}(x))-F(w_{\delta_{0}}(x))\rt)w^1_{\ep,1}(x)\red{\zeta(x)} dx \rt|< \frac{\alpha}{2}.
			\end{equation}
			By \eqref{L6.2}, there exists $\ep_{1}(=\ep_{1}(\alpha))$ such that, for all $\ep<\ep_{1}$, we have
			\begin{equation}\label{L6.7}
			\lt|\int_{\Omega'} (1-\lt|w_{\ep}(x)\rt|^2)F(w_{\delta_{0}}(x))w^1_{\ep,1}(x)\red{\zeta(x)} dx \rt| < \frac{\alpha}{2}.
			\end{equation}
			Define $\ep_{\alpha}:=\min\{\ep_{0},\ep_{1}\}$. Combining (\ref{L6.6}), (\ref{L6.7}), we have that, for all $\ep<\ep_{\alpha}$
			\begin{equation}\label{L6.8}
			\lt|\int_{\Omega'} \lt(1-\lt|w_{\ep}(x)\rt|^2\rt)F(w_{\ep}(x))w^1_{\ep,1}(x)\red{\zeta(x)} dx \rt|< \alpha.
			\end{equation}
			This implies \eqref{L6.3}. The estimates \eqref{L6.4} \blue{and (\ref{L7.4})} follow exactly the same lines.
		\end{proof}
		
		\emph{Proof of Theorem \ref{T4} completed.} Now we return to \eqref{eqs72}. We have
		\begin{equation}\label{T4.7}
		\Pi_{\ep}=\int_{\Omega'} \lt(1-|w_{\ep}|^2\rt)\lt[\Psi_{1,1}(w_{\ep})w_{\ep,1}^1+\Psi_{1,2}(w_{\ep})w_{\ep,1}^2+\Psi_{2,1}(w_{\ep})w_{\ep,2}^1+\Psi_{2,2}(w_{\ep})w_{\ep,2}^2\rt] \red{\zeta(x)}  dx.
		\end{equation}
		\blue{By Lemma \ref{L6} (\ref{L6.3})-(\ref{L6.4}), we have
			\begin{equation}\label{T4.8}
			\int_{\Omega'} \lt(1-|w_{\ep}|^2\rt)\lt[\Psi_{1,1}(w_{\ep})w_{\ep,1}^1+\Psi_{2,2}(w_{\ep})w_{\ep,2}^2\rt] \zeta(x)  dx \rightarrow 0 \text{ as } \ep\rightarrow 0.
			\end{equation}
			By Lemma \ref{l17}, the assumption \eqref{T4.4} implies $\Psi_{1,2}(w_{\ep})=\Psi_{2,1}(w_{\ep})$. Therefore, we have
			\begin{equation}\label{T4.6}
			\int_{\Omega'}\lt(1-|w_{\ep}|^2\rt)\lt[\Psi_{1,2}(w_{\ep})w_{\ep,1}^2+\Psi_{2,1}(w_{\ep})w_{\ep,2}^1\rt] \zeta(x)  dx=\int_{\Omega'}\lt(1-|w_{\ep}|^2\rt)\Psi_{1,2}(w_{\ep})\lt(w_{\ep,1}^2+w_{\ep,2}^1\rt) \zeta(x)  dx.
			\end{equation}
			Now applying Lemma \ref{L6} (\ref{L7.4}) to (\ref{T4.6}) implies that
			\begin{equation}\label{T4.9}
			\int_{\Omega'}\lt(1-|w_{\ep}|^2\rt)\lt[\Psi_{1,2}(w_{\ep})w_{\ep,1}^2+\Psi_{2,1}(w_{\ep})w_{\ep,2}^1\rt] \zeta(x)  dx \rightarrow 0 \text{ as } \ep\rightarrow 0.
			\end{equation}
			Finally, putting (\ref{T4.8}) and (\ref{T4.9}) into (\ref{T4.7}), we obtain $\Pi_{\ep}\rightarrow 0$ as $\ep\rightarrow 0$. This together with (\ref{eqs72}) and (\ref{eqs511}) completes the proof of Theorem \ref{T4}. } 
	\end{proof}

	%%%%%%%%%%%%%%%%%%%%%%%%%%%%%%%%%%%%%%%%%%%%%%%%%%%%%%%
	%
	%   Vanishing of the special entropies
	%
	%%%%%%%%%%%%%%%%%%%%%%%%%%%%%%%%%%%%%%%%%%%%%%%%%%%%%%%%%%%%

	\section{Vanishing of the special entropies}

	Given $\xi\in\mathbb{S}^1$, \bblue{recall the definition of the function $\Phi^{\xi}$ in (\ref{s5.1}). The main result of this section is the following theorem.}
	
	\begin{a2}\label{T5}
		Let $\Omega \subset \R^2$ be a \bblue{bounded simply-connected} domain and \bblue{$u\in E(\Omega)$, where $E(\Omega)$ is defined in (\ref{eqz4})}. Then for every $\xi\in\bblue{\mathbb{S}^1\setminus\{e_1,-e_1,e_2,-e_2\}}$, we have that
		\begin{equation}\label{T5.4}
		\nabla\cdot[\Phi^{\xi}(w)]=0 \qd\text{ in the sense of distributions,}
		\end{equation}
		where $w(x)=(-u_{,2}(x),u_{,1}(x))$.
	\end{a2}
	
	We first recall the following lemma from \cite{mul2}.
	
	\begin{a1}[\cite{mul2}, Lemma 4]\label{l16}
		For a fixed $\xi\in\mathbb{S}^1$, \xblue{the map $\Phi^{\xi}$ defined in (\ref{s5.1})} is a generalized entropy in the sense that there exists a sequence $\{\Phi_{\nu}\}_{\nu \rightarrow \infty}$ of entropies in $C_{c}^{\infty}(\Rb^2;\Rb^2)$ such that
		\begin{equation*} %\label{l16.4}
		\{\Phi_{\nu}(z)\}_{\nu\rightarrow \infty} \text{ is bounded uniformly for bounded } z,
		\end{equation*}
		\begin{equation}\label{l16.5}
		\Phi_{\nu}(z)\rightarrow \Phi^{\xi}(z) \text{ for all } z.
		\end{equation}
	\end{a1}
	
	\xblue{For the convenience of the reader, we include the proof of Lemma \ref{l16} in the Appendix.} Now we provide the proof of Theorem \ref{T5}.
	
	\begin{proof}[Proof of Theorem \ref{T5}]
		Given $\xi\in\bblue{\mathbb{S}^1\setminus\{e_1,-e_1,e_2,-e_2\}}$, we may approximate $\Phi^{\xi}$ by smooth entropies $\Phi^k$ as in Lemma \ref{l16}. We prove that
		\begin{equation*}
			\nabla\cdot[\Phi^k(w)]=0 \text{ in the sense of distributions}
		\end{equation*}
		for all $k$ sufficiently large. As a result, we have \eqref{T5.4}. 
		
		As can be understood from the proof of Theorem \ref{T4}, by virtue of Lemma \ref{L6} the only thing we need to show is
		\begin{equation}\label{t58}
		\int_{\Omega'}\lt(1-|w_{\ep}|^2\rt)\lt[\Psi^k_{1,2}(w_{\ep})w_{\ep,1}^2+\Psi^k_{2,1}(w_{\ep})w_{\ep,2}^1\rt]\zeta dx \rightarrow 0,
		\end{equation}
		where the function $\Psi^k$ is related to $\Phi^k$ through Lemma \ref{l12} and $\zeta\in C^{\infty}_{c}(\Omega')$ is any test function. Let us write
		\begin{equation*}
		\begin{split}
		&\int_{\Omega'}\lt(1-|w_{\ep}|^2\rt)\lt[\Psi^k_{1,2}(w_{\ep})w_{\ep,1}^2+\Psi^k_{2,1}(w_{\ep})w_{\ep,2}^1\rt]\zeta dx\\
		=&\int_{\Omega'}\lt(1-|w_{\ep}|^2\rt)\frac{\Psi^k_{1,2}(w_{\ep})+\Psi^k_{2,1}(w_{\ep})}{2}\lt(w_{\ep,1}^2+w_{\ep,2}^1\rt)\zeta dx\\
		&\quad\quad\quad\quad\quad\quad\quad\quad\quad\quad\quad\quad+\int_{\Omega'}\lt(1-|w_{\ep}|^2\rt)\frac{\Psi^k_{1,2}(w_{\ep})-\Psi^k_{2,1}(w_{\ep})}{2}\lt(w_{\ep,1}^2-w_{\ep,2}^1\rt)\zeta dx.
		\end{split}
		\end{equation*}
		We deduce from Lemma \ref{L6} (\ref{L7.4}) that
			\begin{equation}
			\label{t51}
			\int_{\Omega'}\lt(1-|w_{\ep}|^2\rt)\lt(w_{\ep,1}^2+w_{\ep,2}^1\rt)\frac{\Psi^k_{1,2}(w_{\ep})+\Psi^k_{2,1}(w_{\ep})}{2}\zeta dx \rightarrow 0 \quad \text{ as } \ep\rightarrow 0.
			\end{equation}
		In the following, we show
		\begin{equation}\label{t52}
		\int_{\Omega'}\lt(1-|w_{\ep}|^2\rt)\frac{\Psi^k_{1,2}(w_{\ep})-\Psi^k_{2,1}(w_{\ep})}{2}\lt(w_{\ep,1}^2-w_{\ep,2}^1\rt)\zeta dx\rightarrow 0 \quad \text{ as } \ep\rightarrow 0.
		\end{equation}
		
		Let us denote $\psi^k(z)=\frac{\nabla\lt(\Delta\varphi^{k}\rt)\cdot z^{\perp}}{4}$, where the function $\varphi^{k}$ is related to $\Phi^{k}$ through Lemma \ref{l13}. Using this new function $\psi^k$ and the calculation \eqref{s4.2}, we write
		\begin{equation}\label{t53}
		\begin{split}
		&\int_{\Omega'}\lt(1-|w_{\ep}|^2\rt)\frac{\Psi^k_{1,2}(w_{\ep})-\Psi^k_{2,1}(w_{\ep})}{2}\lt(w_{\ep,1}^2-w_{\ep,2}^1\rt)\zeta dx\\
		=&\int_{\Omega'}\lt(1-|w_{\ep}|^2\rt)\psi^k(w)\lt(w_{\ep,1}^2-w_{\ep,2}^1\rt)\zeta dx + \int_{\Omega'}\lt(1-|w_{\ep}|^2\rt)\lt(\psi^k(w_{\ep})-\psi^k(w)\rt)\lt(w_{\ep,1}^2-w_{\ep,2}^1\rt)\zeta dx.
		\end{split}
		\end{equation}
		Recall the definition of $w_{\ep}$ in \eqref{s4.4}. For the above first term, we further write
		\begin{equation}\label{t54}
		\begin{split}
		&\int_{\Omega'}\lt(1-|w_{\ep}|^2\rt)\psi^k(w)\lt(w_{\ep,1}^2-w_{\ep,2}^1\rt)\zeta dx\\
		=&\int_{\Omega'}\lt(1-|\nabla u_{\ep}|^2\rt)\psi^k(w)\lt(u_{\ep,11}+u_{\ep,22}\rt)\zeta dx\\
		=&\int_{\Omega'}\lt(1-|\nabla u_{\ep}|^2\rt)\psi^k(w)\lt(u_{\ep,11}+u_{\ep,22}+\frac{u_{\ep,12}}{u_{,1}u_{,2}}\rt)\zeta dx-\int_{\Omega'}\lt(1-|\nabla u_{\ep}|^2\rt)u_{\ep,12}\frac{\psi^k(w)}{u_{,1}u_{,2}}\zeta dx.
		\end{split}
		\end{equation}
		In the following, we will establish 
		\begin{equation}\label{t55}
		\int_{\Omega'}\lt(1-|w_{\ep}|^2\rt)\lt(\psi^k(w_{\ep})-\psi^k(w)\rt)\lt(w_{\ep,1}^2-w_{\ep,2}^1\rt)\zeta dx\rightarrow 0,
		\end{equation}
		\begin{equation}\label{t56}
		\int_{\Omega'}\lt(1-|\nabla u_{\ep}|^2\rt)\psi^k(w)\lt(u_{\ep,11}+u_{\ep,22}+\frac{u_{\ep,12}}{u_{,1}u_{,2}}\rt)\zeta dx\rightarrow 0,
		\end{equation}
		and
		\begin{equation}\label{t57}
		\int_{\Omega'}\lt(1-|\nabla u_{\ep}|^2\rt)u_{\ep,12}\frac{\psi^k(w)}{u_{,1}u_{,2}}\zeta dx\rightarrow 0
		\end{equation}
		as $\ep\rightarrow 0$, respectively. Putting \eqref{t53}-\eqref{t57} together, we obtain \eqref{t52}, which together with \eqref{t51} gives us \eqref{t58}. This will conclude the proof of the theorem.
		
		\red{First note (\ref{t55}) follows as a direct consequence of Lemma \ref{l14}. Equations (\ref{t56}), (\ref{t57}) will be established in the following two lemmas. } 
	\end{proof}
	
	\red{\begin{a1}
			\label{l10}
			We have 
			\begin{equation}
			\label{t57.5}
			\int_{\Omega'}\lt(1-|\nabla u_{\ep}|^2\rt)u_{\ep,12}\frac{\psi^k(w)}{u_{,1}u_{,2}}\zeta dx\rightarrow 0 \blue{\text{ as } \ep\rightarrow 0.}
			\end{equation}
		\end{a1}}
		
		\begin{proof}
			A key observation in the proof is that, for $k$ sufficiently large, $\chi^k:=\frac{\psi^k(w)}{u_{,1}u_{,2}}$ is an $L^{\infty}$ function. Indeed, we use smooth entropies $\Phi^k$ to approximate the entropy $\Phi^{\xi}$ in the way that is given in the proof of Lemma \ref{l16} in \xblue{the Appendix}. In particular, for $k$ sufficiently large, the function $\varphi^k$ satisfies \red{$D^2 \varphi^k=0$} outside a sufficiently small neighborhood of the line \xblue{$z\cdot \xi =0$} inside the ball $B_k(0)$. Consequently, \emph{on $\mathbb{S}^1$}, $\psi^k(z)=\frac{\nabla\lt(\Delta\varphi^{k}\rt)\cdot z^{\perp}}{4}$ is supported in a sufficiently small neighborhood of \xred{ the points }\xblue{$z\cdot \xi =0$ with $|z|=1$}. Since we have chosen $\xi\in\mathbb{S}^1$ to be such that $\xi$ is not parallel to the axes, \xblue{for $k$ sufficiently large,} the support of $\psi^k(z)$ \emph{on $\mathbb{S}^1$} \xblue{is bounded away from} the axes. Indeed, let $\alpha>0$ denote the distance between the support of $\psi^{k}$ on $\mathbb{S}^1$ and the axes. Then, either $|u_{,i}|<\frac{\alpha}{2}$ for some $i=1,2$, so $\psi^{k}(w)=0$, or $|u_{,1}|\geq \frac{\alpha}{2}$ and $|u_{,2}|\geq \frac{\alpha}{2}$, so $|\chi^k|\leq \frac{4\lVert\psi^{k}\rVert_{\infty}}{\alpha^2}$. Therefore, for all $x\in\Omega$ such that $|\nabla u(x)|=1$, we have $\chi^k(x)\leq C_{k}$ for some constant $C_{k}$ depending only on $\Phi^k$. Since $|\na u|=1$ a.e. in $\Omega$, we have $\chi^k\in L^{\infty}(\Omega)$.
			
			In particular, we have $\chi^k\in L^{4}(\Omega)$. Let $\{\chi_j\}$ be a sequence of smooth functions such that
			\begin{equation*} %\label{l101}
			\chi_j\rightarrow \chi^k \quad\text{ in } L^{4}(\Omega).
			\end{equation*}
			Then we have
			\begin{equation}\label{l102}
			\begin{split}
			&\int_{\Omega'}\lt(1-|\nabla u_{\ep}|^2\rt)u_{\ep,12}\frac{\psi^k(w)}{u_{,1}u_{,2}}\zeta dx\\
			=&\int_{\Omega'}\lt(1-|\nabla u_{\ep}|^2\rt)u_{\ep,12}\chi_j\zeta dx+\int_{\Omega'}\lt(1-|\nabla u_{\ep}|^2\rt)u_{\ep,12}\lt(\chi^k-\chi_j\rt)\zeta dx.
			\end{split}
			\end{equation}
			It follows from Lemma \ref{l11} that
			\begin{equation}\label{l103}
			\int_{\Omega'}\lt(1-|\nabla u_{\ep}|^2\rt)u_{\ep,12}\lt(\chi^k-\chi_j\rt)\zeta dx\rightarrow 0 \quad \text{ as } \ep\rightarrow 0, j\rightarrow\infty.
			\end{equation}
			On the other hand, we have $\chi_j\zeta\in C^{\infty}_{c}(\Omega')$. It follows from Lemma \ref{l15} (noting the relationship between $w_{\ep}$ and $u_{\ep}$ as in \eqref{s4.4}) that
			\begin{equation}\label{l104}
			\int_{\Omega'}\lt(1-|\nabla u_{\ep}|^2\rt)u_{\ep,12}\chi_j\zeta dx\rightarrow 0 \quad \text{ as } \ep\rightarrow 0.
			\end{equation}
			Putting \eqref{l102}-\eqref{l104} together \blue{ and using arguments similar to those in (\ref{L6.6})-(\ref{L6.8})}, we obtain \eqref{t57.5}.
		\end{proof}
		
		\begin{a1}\label{l9}
			We have	
			\red{\begin{equation*}
				%\label{t56.5}
				\int_{\Omega'}\lt(1-|\nabla u_{\ep}|^2\rt)\psi^k(w)\lt(u_{\ep,11}+u_{\ep,22}+\frac{u_{\ep,12}}{u_{,1}u_{,2}}\rt)\zeta dx\rightarrow 0 \blue{\text{ as } \ep\rightarrow 0.}
				\end{equation*}}
		\end{a1}
		
		\begin{proof}
			Recall that we defined $\chi^k:=\frac{\psi^k(w)}{u_{,1}u_{,2}}\in L^{\infty}(\Omega)$. We write
			\begin{equation}\label{l96}
			\begin{split}
			&\int_{\Omega'}\lt(1-|\nabla u_{\ep}|^2\rt)\psi^k(w)\lt(u_{\ep,11}+u_{\ep,22}+\frac{u_{\ep,12}}{u_{,1}u_{,2}}\rt)\zeta dx\\
			=&\int_{\Omega'}\lt(1-|\nabla u_{\ep}|^2\rt)\lt(u_{,1}u_{,2}\lt(u_{\ep,11}+u_{\ep,22}\rt)+u_{\ep,12}\rt)\overbrace{\frac{\psi^k(w)}{u_{,1}u_{,2}}}^{\chi^k}\zeta dx\\
			=&\overbrace{\int_{\Omega'}\lt(1-|\nabla u_{\ep}|^2\rt)\lt(u_{\ep,1}u_{\ep,2}\lt(u_{\ep,11}+u_{\ep,22}\rt)+u_{\ep,12}\rt)\chi^k\zeta dx}^{I}+\overbrace{\int_{\Omega'}\lt(1-|\nabla u_{\ep}|^2\rt)\lt(u_{\ep,11}+u_{\ep,22}\rt)\lt(u_{,1}u_{,2}-u_{\ep,1}u_{\ep,2}\rt)\chi^k\zeta dx}^{II}\\
			=&\overbrace{\int_{\Omega'}\lt(1-|\nabla u_{\ep}|^2\rt)\lt(u_{\ep,1}u_{\ep,2}\lt(u_{\ep,11}+u_{\ep,22}\rt)+|\nabla u_{\ep}|^2u_{\ep,12}\rt)\chi^k\zeta dx+\int_{\Omega'}\lt(1-|\nabla u_{\ep}|^2\rt)^2u_{\ep,12}\chi^k\zeta dx}^{I}\\
			&+\overbrace{\int_{\Omega'}\lt(1-|\nabla u_{\ep}|^2\rt)\lt(u_{\ep,11}+u_{\ep,22}\rt)\lt(u_{,1}u_{,2}-u_{\ep,1}u_{\ep,2}\rt)\chi^k\zeta dx}^{II}.
			\end{split}
			\end{equation}
			First, we have (noting $\lt|\nabla u\rt|=1$ a.e.)
			\begin{equation}\label{l91}
			\int_{\Omega'}\lt(1-|\nabla u_{\ep}|^2\rt)^2u_{\ep,12}\chi^k\zeta dx=\int_{\Omega'}\lt(|\nabla u|^2-|\nabla u_{\ep}|^2\rt)\lt(1-|\nabla u_{\ep}|^2\rt)u_{\ep,12}\chi^k\zeta dx.
			\end{equation}
			Since $|\nabla u_{\ep}|\leq1$ and $|\nabla u|=1$, and $\chi^k \zeta\in \blue{L^{\infty}(\Omega)}$, we have
			\begin{equation}\label{l92}
			\lt|\int_{\Omega'}\lt(|\nabla u|^2-|\nabla u_{\ep}|^2\rt)\lt(1-|\nabla u_{\ep}|^2\rt)u_{\ep,12}\chi^k\zeta dx\rt| \leq C\int_{\Omega'}\lt|\nabla u-\nabla u_{\ep}\rt|\lt(1-|\nabla u_{\ep}|^2\rt)\lt|u_{\ep,12}\rt| dx.
			\end{equation}
			Since $\lVert \nabla u-\nabla u_{\ep} \rVert_{L^{4}(\Omega')}=\lt\lVert w-w_{\ep}\rt\rVert_{L^{4}(\Omega')} \rightarrow 0$, we deduce from \eqref{l92} and Lemma \ref{l11} that
			\begin{equation}\label{l93}
			\int_{\Omega'}\lt(|\nabla u|^2-|\nabla u_{\ep}|^2\rt)\lt(1-|\nabla u_{\ep}|^2\rt)u_{\ep,12}\chi^k\zeta dx\rightarrow 0 \quad \text{ as } \ep \rightarrow 0.
			\end{equation}
			Combining \eqref{l91} with \eqref{l93}, we obtain
			\begin{equation}\label{l94}
			\int_{\Omega'}\lt(1-|\nabla u_{\ep}|^2\rt)^2u_{\ep,12}\chi^k\zeta dx\rightarrow 0 \quad \text{ as } \ep \rightarrow 0.
			\end{equation}
			
			For the last term in \eqref{l96}, \blue{since $\lVert w-w_{\ep}\rVert_{L^{p}(\Omega')}\rightarrow 0$ for all $p\geq 1$, it is clear that $\lVert u_{,1}u_{,2}-u_{\ep,1}u_{\ep,2}\rVert_{L^{4}(\Omega')}\rightarrow 0$}. It follows from \blue{the fact that $\chi^k \zeta\in L^{\infty}(\Omega)$ and} Lemma \ref{l11} again that
			\begin{equation}\label{l95}
			\int_{\Omega'}\lt(1-|\nabla u_{\ep}|^2\rt)\lt(u_{\ep,11}+u_{\ep,22}\rt)\lt(u_{,1}u_{,2}-u_{\ep,1}u_{\ep,2}\rt)\chi^k\zeta dx\rightarrow 0 \quad \text{ as } \ep \rightarrow 0.
			\end{equation}
			
			Finally, we look at the first term in \eqref{l96}. Following the arguments in Lemma \ref{l10}, we choose a sequence of smooth functions $\{\chi_j\}$ such that
			\begin{equation*}
				\chi_j\rightarrow \chi^k \quad\text{ in } L^{4}(\Omega).
			\end{equation*} 
			Note that we have $|u_{\ep,1}u_{\ep,2}|\leq 1$ and $|\nabla u_{\ep}|\leq 1$. Therefore, we have
			\begin{equation*}
				\begin{split}
					&\lt|\int_{\Omega'}\lt(1-|\nabla u_{\ep}|^2\rt)\lt(u_{\ep,1}u_{\ep,2}\lt(u_{\ep,11}+u_{\ep,22}\rt)+|\nabla u_{\ep}|^2u_{\ep,12}\rt)\lt(\chi^k-\chi_j\rt)\zeta dx\rt|\\
					&\quad\quad\quad\quad\quad\quad\quad\quad\quad\quad \leq\int_{\Omega'}\lt(1-|\nabla u_{\ep}|^2\rt)\lt(\lt|u_{\ep,11}\rt|+\lt|u_{\ep,22}\rt|+\lt|u_{\ep,12}\rt|\rt)\lt|\chi^k-\chi_j\rt|\lt|\zeta\rt| dx.
				\end{split}
			\end{equation*}
			\red{ Let $\alpha>0$. By Lemma \ref{l11} there exists some $j_0\in \mathbb{N}$ such that 
				\begin{equation}\label{l97}
				\lt|\int_{\Omega'}\lt(1-|\nabla u_{\ep}|^2\rt)\lt(u_{\ep,1}u_{\ep,2}\lt(u_{\ep,11}+u_{\ep,22}\rt)+|\nabla u_{\ep}|^2u_{\ep,12}\rt)\lt(\chi^k-\chi_j\rt)\zeta dx\rt|\leq \frac{\alpha}{2}\text{ for all }\ep\in (0,\ep_0), j\geq j_0
				\end{equation}
				where $\ep_{0}$ is the small constant as in Lemma \ref{l11}}.
			
			Using the harmonic polynomial $\xblue{\wt\varphi(z)}=z_1^2-z_2^2$ and the formula \eqref{s4.1}, we obtain $\xblue{\wt{\Phi}(z)}=(z_1^3+3z_1z_2^2,-3z_1^2z_2-z_2^3)$. Let $\green{\eta}\in C^{\infty}_{c}(\Rb^2)$ be a cut-off function such that $\green{\eta}\equiv 1$ on $B_2(0)$ and define $\varphi:=\wt\varphi\green{\eta}\in C^{\infty}_{c}(\R^2)$. Let $\Phi$ be the entropy obtained from the function $\varphi$ through formula (\ref{s4.1}). Since $\varphi=\wt\varphi$ on $B_{2}(0)$, we have $\Phi=\wt\Phi$ on $B_{2}(0)$. Since $|w|=1$ a.e. \xblue{and $|w_{\ep}|\leq 1$}, we see that \xblue{$\Phi(w)=\wt\Phi(w)$ for a.e. $x\in\Omega$ and $\Phi(w_{\ep})=\wt\Phi(w_{\ep})$ for all $x\in\Omega'$.} By direct calculations, we have
			\begin{equation*}
				\begin{split}
					\nabla\cdot[\Phi(w_{\ep})]&=\nabla\cdot\lt(-u_{\ep,2}^3-3u_{\ep,2}u_{\ep,1}^2,-3u_{\ep,2}^2u_{\ep,1}-u_{\ep,1}^3\rt)\\
					&=-6\lt(u_{\ep,1}u_{\ep,2}\lt(u_{\ep,11}+u_{\ep,22}\rt)+|\nabla u_{\ep}|^2u_{\ep,12}\rt).
				\end{split}
			\end{equation*} 
			Let us apply \eqref{eqs72} to our particular entropy $\Phi$:
			\begin{equation}\label{l912}
			\begin{split}
			&\int_{\Omega'}\lt(1-|\nabla u_{\ep}|^2\rt)\lt(u_{\ep,1}u_{\ep,2}\lt(u_{\ep,11}+u_{\ep,22}\rt)+|\nabla u_{\ep}|^2u_{\ep,12}\rt)\chi_{j_0}\zeta dx\\
			=&-\frac{1}{6}\int_{\Omega'}\lt(1-|\nabla u_{\ep}|^2\rt)\chi_{j_0}\zeta\nabla\cdot\lt[\Phi(w_{\ep})\rt] dx\\
			\overset{(\ref{l18.1})}{=}&-\frac{1}{6}\int_{\Omega'}\lt(1-|\nabla u_{\ep}|^2\rt)\chi_{j_0}\zeta\Psi(w_{\ep})\cdot\nabla\lt(1-|w_{\ep}|^2\rt) dx\\
			=&-\frac{1}{12}\int_{\Omega'}\chi_{j_0}\zeta\Psi(w_{\ep})\cdot\nabla\lt(1-|w_{\ep}|^2\rt)^2 dx\\
			=&-\frac{1}{12}\int_{\Omega'} \chi_{j_0}\zeta \na \cdot \lt[\Psi(w_{\ep}) \lt(1-\lt|w_{\ep}\rt|^2 \rt)^2  \rt] dx+\frac{1}{12}\int_{\Omega'}\chi_{j_0}\zeta  \lt(1-\lt|w_{\ep}\rt|^2 \rt)^2 \na\cdot \lt[ \Psi(w_{\ep}) \rt] dx,
			\end{split}
			\end{equation}
			\xblue{where $\Psi \in C_{c}^{\infty}(\Rb^2;\Rb^2)$ is related to the particular entropy $\Phi$ via Lemma \ref{l12}.} It is clear that \xblue{$\{\sup|\Psi(w_{\ep})|\}$} is uniformly bounded. It follows from integration by parts that
			\begin{equation}\label{l913}
			\int_{\Omega'} \chi_{j_0}\zeta \na \cdot \lt[\Psi(w_{\ep}) \lt(1-\lt|w_{\ep}\rt|^2 \rt)^2  \rt] dx\rightarrow 0 \quad \text{ as } \ep \rightarrow 0.
			\end{equation}
			Now we write out the other term in (\ref{l912})
			\begin{equation}\label{l98}
			\begin{split}
			&\int_{\Omega'}\chi_{j_0}\zeta  \lt(1-\lt|w_{\ep}\rt|^2 \rt)^2 \na\cdot \lt[ \Psi(w_{\ep}) \rt] dx\\
			=&\int_{\Omega'}\chi_{j_0}\zeta  \lt(1-\lt|w_{\ep}\rt|^2 \rt)^2\lt[\Psi_{1,1}(w_{\ep})w_{\ep,1}^1+\Psi_{1,2}(w_{\ep})w_{\ep,1}^2+\Psi_{2,1}(w_{\ep})w_{\ep,2}^1+\Psi_{2,2}(w_{\ep})w_{\ep,2}^2\rt] dx.
			\end{split}
			\end{equation}
			\xblue{For all $m,n\in\{1,2\}$}, using $|w|=1$ a.e., we have
			\begin{equation}\label{l99}
			\begin{split}
			\int_{\Omega'}\chi_{j_0}\zeta  \lt(1-\lt|w_{\ep}\rt|^2 \rt)^2\Psi_{\xblue{m,n}}(w_{\ep})w_{\ep,\xblue{m}}^{\xblue{n}}dx=\int_{\Omega'}\chi_{j_0}\zeta  \lt(\lt|w\rt|^2-\lt|w_{\ep}\rt|^2 \rt)\lt(1-\lt|w_{\ep}\rt|^2 \rt)\Psi_{m,n}(w_{\ep})w_{\ep,m}^ndx.
			\end{split}
			\end{equation}
			Since $\lVert w-w_{\ep} \rVert_{L^{4}(\Omega')}\rightarrow 0$ and \xblue{$\{\sup|\Psi_{m,n}(w_{\ep})|\}$} is uniformly bounded, an application of Lemma \ref{l11} yields
			\begin{equation}\label{l911}
			\int_{\Omega'}\chi_{j_0}\zeta  \lt(\lt|w\rt|^2-\lt|w_{\ep}\rt|^2 \rt)\lt(1-\lt|w_{\ep}\rt|^2 \rt)\Psi_{m,n}(w_{\ep})w_{\ep,m}^ndx\leq C\int_{\Omega'}\lt|w-w_{\ep}\rt|\lt(1-\lt|w_{\ep}\rt|^2 \rt)w_{\ep,m}^ndx\rightarrow 0.
			\end{equation}
			Putting \eqref{l99}-\eqref{l911} together, we obtain
			\begin{equation*}
				\int_{\Omega'}\chi_{j_0}\zeta  \lt(1-\lt|w_{\ep}\rt|^2 \rt)^2\Psi_{m,n}(w_{\ep})w_{\ep,m}^ndx\rightarrow 0.
			\end{equation*}
			\xblue{Taking the sum over all $m, n$,} we deduce from \eqref{l98} that
			\begin{equation}\label{l914}
			\int_{\Omega'}\chi_{j_0}\zeta  \lt(1-\lt|w_{\ep}\rt|^2 \rt)^2 \na\cdot \lt[ \Psi(w_{\ep}) \rt] dx\rightarrow 0.
			\end{equation}
			Combining \eqref{l914} with \eqref{l913} and \eqref{l912}, we have
			\begin{equation*}
				\int_{\Omega'}\lt(1-|\nabla u_{\ep}|^2\rt)\lt(u_{\ep,1}u_{\ep,2}\lt(u_{\ep,11}+u_{\ep,22}\rt)+|\nabla u_{\ep}|^2u_{\ep,12}\rt)\chi_{j_0}\zeta dx\rightarrow 0\text{ as }\ep\rightarrow 0.
			\end{equation*}
			
			\red{So there exists some $\ep_1\in (0,\ep_0)$ such that 
				\begin{equation}
				\label{eqzz91}
				\int_{\Omega'}\lt(1-|\nabla u_{\ep}|^2\rt)\lt(u_{\ep,1}u_{\ep,2}\lt(u_{\ep,11}+u_{\ep,22}\rt)+|\nabla u_{\ep}|^2u_{\ep,12}\rt)\chi_{j_0}\zeta dx<\frac{\alpha}{2}\text{ for any }\ep\in \lt(0,\ep_1\rt).
				\end{equation}
				Inequality (\ref{eqzz91}) together with \eqref{l97} yield
				\begin{equation}\label{l915.5}
				\lt|\int_{\Omega'}\lt(1-|\nabla u_{\ep}|^2\rt)\lt(u_{\ep,1}u_{\ep,2}\lt(u_{\ep,11}+u_{\ep,22}\rt)+|\nabla u_{\ep}|^2u_{\ep,12}\rt)\chi^k\zeta dx\rt|<\alpha \text{ for any }\ep\in (0,\ep_1).\nn
				\end{equation}
				As this is true for any $\alpha>0$ we have shown
				\begin{equation}\label{l915}
				\int_{\Omega'}\lt(1-|\nabla u_{\ep}|^2\rt)\lt(u_{\ep,1}u_{\ep,2}\lt(u_{\ep,11}+u_{\ep,22}\rt)+|\nabla u_{\ep}|^2u_{\ep,12}\rt)\chi^k\zeta dx\rightarrow 0\text{ as }\ep\rightarrow 0. 
				\end{equation}}
			Finally, putting \eqref{l94}, \eqref{l95} and \eqref{l915} into \eqref{l96} concludes the proof of Lemma \ref{l9}.
		\end{proof}
		
%%%%%%%%%%%%%%%%%%%%%%%%%%%%%%%%%%%%%%%%%%%%%%%%%%%%%%%%%%%%%%%%%%%%%%%%%%%%%%%%%%%%%%%%%%%%%%%
%
%
%%%%%%%%%%%%%%%%%%%%%%%%%%%%%%%%%%%%%%%%%%%%%%%%%%%%%%%%%%%%%%%%%%%%%%%%%%%%%%%%%%%%%%%%%%%%%%%		
		
\section{Proof of Theorem \ref{T1}}

By Theorem \ref{T5}, (\ref{T5.4}) we have that 
\begin{equation}
	\label{eqj11}
	\nabla\cdot\lt[\Phi^{\xi}\lt(\na u^{\perp}\rt)\rt]=0\text{ distributionally in }\Omega,\text{ for any }\xi\in \mathbb{S}^1\backslash \lt\{e_1,-e_1,e_2,-e_2\rt\}.
\end{equation} 
As explained in the sketch of the proof, we could carry out the the argument that establishes \green{(\ref{eqj11})} for a coordinate axis $\lt\{\ep_1,\ep_2\rt\}$ (see (\ref{eqzzz222})) and this gives (\ref{eqj11}) for all $\xi\in \mathbb{S}^1\backslash \lt\{\ep_1,-\ep_1,\ep_2,-\ep_2\rt\}$ and hence (\ref{eqj11}) holds for 
any $\xi\in \mathbb{S}^1$. 

Now defining $\bblue{w}(x)=\na u(x)^{\perp}$ we have that $\Phi^{\xi}\lt(\na u(x)^{\perp}\rt)\overset{(\ref{eqi45.6})}{=} \xi\chi\lt(x,\xi\rt)$ for a.e.\ $x\in \Omega$ and so 
$$
0=\nabla\cdot\lt[\Phi^{\xi}\lt(\na u^{\perp}\rt)\rt]=\nabla\cdot\lt[\xi \chi\lt(\cdot,\xi\rt)\rt]=\xi\cdot \na \chi\lt(\cdot,\xi\rt)\text{ in }\mathcal{D}'(\Omega)
$$
and thus applying Theorem \ref{JOP1} we have that $\na u$ is \bblue{locally} Lipschitz outside a \bblue{locally} finite set of points. 
		
It has been observed in \cite{ignat} that the results of \cite{otto} imply that under the hypothesis of Theorem \ref{JOP1}, if $\OI\subset \subset \Omega$ is a convex neighborhood of a point $\zeta\in S$ (where $\bblue{w}\green{=\na u^{\perp}}$ is locally Lipschitz outside of $S$) then \green{there exists $\alpha\in\{1,-1\}$ such that}
\begin{equation}
\label{eqj71}
\bblue{w}(z)=\green{\alpha}\frac{(z-\zeta)}{\lt|z-\zeta\rt|}^{\perp}\text{ for any }z\in \OI.
\end{equation}
Since we have shown that $\bblue{w}$ satisfies (\ref{eqi45}), this implies (\ref{eq210}). For the convenience of the reader, we note that (\ref{eqj71}) follows from the results of \cite{otto} in the following way. Firstly by Lemma 5.1 \cite{otto} for any $x_0,y_0\in \OI$ that are Lebegue points of $\bblue{w}$ we have 
\begin{equation}
\label{eqj72}
\lt|\bblue{w}(x_0)-\alpha \bblue{w}(y_0)\rt|\leq \frac{\lt|x_0-y_0\rt|}{d}\text{ for some }\alpha\in \lt\{1,-1\rt\},
\end{equation}
\xblue{where $d=\mathrm{dist}(\OI,\partial\Omega)>0$.} In the proof of Theorem 1.3 (that follows the proof of Lemma 5.1) the estimate (\ref{eqj72}) is strengthened in that it is shown that $\alpha=1$. Thus $\bblue{w}$ is $\frac{1}{d}$-Lipschitz in $\OI$. This contradicts the fact that $\zeta\in \OI$ and hence (\ref{eqj71}) follows. 	$\Box$

		\section{Appendix: Some auxiliary results}
				
				We \green{have used} in a fundamental way a couple of estimates from \cite{DeI}, these in turns were inspired \bblue{by} a commutator estimate of Constantin, E, Titi \cite{titi}. For convenience of the reader we repeat the proof from \cite{DeI}.
				
				\begin{a1}[\cite{DeI}]
				\label{l21}
					Let $\Omega \subset \Rb^2$ be a bounded domain and $w \in L^3(\Omega;\Rb^2)$ satisfy $|w|=1$ a.e. in $\Omega$. \xblue{Given $\Omega'\subset\subset\Omega$, let $\gamma:=\mathrm{dist}(\Omega',\partial\Omega)>0$.} Then, \xblue{for all $x\in\Omega'$ and $0<\ep<\gamma$}, denoting $w_{\ep}=w\ast \rho_{\ep}$, we have
					\begin{equation}\label{l21.1}
					1-\lt|w_{\ep}(x)\rt|^2 \leq \frac{2\lVert \rho \rVert_{L^{\infty}}}{\ep^2}\int_{B_{\ep}}\lt|w(x-z)-w(x)\rt|^2 dz,
					\end{equation}
					and 
					\begin{equation}\label{l21.2}
					|\partial_{j}w_{\ep}(x)| \leq \frac{\lVert \nabla\rho \rVert_{L^{\infty}}}{\ep^3}\int_{B_{\ep}}\lt| w(x-z)-w(x)\rt|dz.
					\end{equation}
				\end{a1}
				
				\begin{proof}
					First, for $x \in \xblue{\Omega'}$ and for \xblue{$0<\ep<\gamma$}, using $|w|=1$ a.e., we have
					\begin{equation*}
						\begin{split}
							1-\lt|w_{\ep}(x)\rt|^2 &= |w|^2\ast\rho_{\ep}(x)-\lt|w\ast\rho_{\ep}\rt|^2\\
							&= \int_{\Rb^2}\lt|w(x-z)\rt|^2 \rho_{\ep}(z)dz\\
							&\qd\qd -\lt(\int_{\Rb^2}w(x-z)\rho_{\ep}(z)dz\rt)\cdot \lt(\int_{\Rb^2}\xblue{w}(x-y)\rho_{\ep}(y)dy\rt)\\
							&= \int_{\Rb^2}\int_{\Rb^2} w(x-z)\lt(w(x-z)-w(x-y)\rt)\rho_{\ep}(z)\rho_{\ep}(y)dz dy\\
							&\overset{z:=y, y:=z}{=}\frac{1}{2}\int_{\Rb^2}\int_{\Rb^2}\lt|w(x-z)-w(x-y)\rt|^2 \rho_{\ep}(z)\rho_{\ep}(y)dz dy\\
							&\leq 2\int_{\Rb^2} \lt|w(x-z)-w(x)\rt|^2 \rho_{\ep}(z) dz\\
							&\leq \frac{2\lVert \rho \rVert_{L^{\infty}}}{\ep^2}\int_{B_{\ep}}\lt|w(x-z)-w(x)\rt|^2 dz.
						\end{split}
					\end{equation*}
					This establishes (\ref{l21.1}).
					
					To show \eqref{l21.2}, note that $\int_{B_{\ep}}\partial_j \rho(\frac{z}{\ep})dz = 0$ for $j=1, 2$. Therefore, we have
					\begin{equation*}
						\begin{split}
							\lt|\partial_j w_{\ep}(x)\rt| &= \lt|w\ast\partial_j \rho_{\ep}(x)\rt| = \lt|\frac{1}{\ep^3} \int_{B_{\ep}}w(x-z)\partial_j\rho(\frac{z}{\ep})dz\rt|\\
							&= \lt|\frac{1}{\ep^3} \int_{B_{\ep}}\lt(w(x-z)-w(x)\rt)\partial_j\rho(\frac{z}{\ep})dz\rt|\\
							&\leq \frac{\lVert \nabla\rho \rVert_{L^{\infty}}}{\ep^3}\int_{B_{\ep}}\lt| w(x-z)-w(x)\rt|dz.
						\end{split}
					\end{equation*}
				\end{proof}

		\begin{a1}
			\label{curl}
			Let $\Omega\subset\R^2$ be a bounded simply-connected domain and $v\in L^{\infty}(\Omega;\R^2)$ be such that $\mathrm{curl}v=0$ weakly. Then there exists some potential $f\in W^{1,\infty}(\Omega)$ such that $\na f=v$ a.e. on $\Omega$. 
		\end{a1}
		
		\begin{proof}[Proof of Lemma \ref{curl}]
		We follow some of the ideas in the proof of Theorem 2.9 in \cite{gr}. The proof goes in two steps.
			
		\em Step 1. \rm We can find a sequence $\{\Omega_k\}_k$ of open simply-connected sets with the following properties:
		\begin{enumerate}
		\item $\Omega_k\subset\subset \Omega$;
		\item $\Omega_k\subset \Omega_{k+1}$;
		\item $\bigcup_{k} \Omega_k=\Omega$.
		\end{enumerate}
		
		\em Proof of Step 1. \rm Define $O_k:=\{x\in\Omega: \mathrm{dist}(x,\partial\Omega)>2^{-k} \}$. We start with some $k_0$ sufficiently large such that $O_{k_0}$ is nonempty. Define $\Omega_{k_0}$ to be any connected component of $O_{k_0}$. For all $k>k_0$, define $\Omega_k$ to be the connected component of $O_k$ that contains $\Omega_{k_0}$. It is clear that the sequence $\{\Omega_k\}_{k\geq k_0}$ satisfies (1) and (2). To see (3), we claim that for any $k_1\geq k_0$, $O_{k_1}\subset \Omega_k$ for all $k$ sufficiently large. Indeed, let $\{O_{k_1}^j\}_{j=1}^{m}$ be the connected components of $O_{k_1}$. Without loss of generality, assume $O_{k_1}^1=\Omega_{k_1}$. For each $j=1,2,...,m$, we fix a point $a_j\in O_{k_1}^j$. Since $\Omega$ is connected, we can find continuous paths $\gamma_1^j\subset\Omega$ connecting $a_1$ and $a_j$ for $j=2,3,...,m$. Denote $\delta_j=\mathrm{dist}(\gamma_1^j,\partial\Omega)>0$, and let $T_j$ be a tubular neighborhood of $\gamma_1^j$ of size $\frac{\delta_j}{2}$ for $j=2,3,...,m$. Now denote $\delta_{k_1}:=\min\{\delta_j\}>0$. It is clear that $O_{k_1}\bigcup\lt(\bigcup_{j=2}^{m}T_j\rt)\subset\Omega$ is connected, and for any $x\in O_{k_1}\bigcup\lt(\bigcup_{j=2}^{m}T_j\rt)$, $\mathrm{dist}(x,\partial\Omega)>\min\{2^{-k_1},\frac{\delta_{k_1}}{2}\}$. Therefore $O_{k_1}\bigcup\lt(\bigcup_{j=2}^{m}T_j\rt)\subset O_k$ for $k$ sufficiently large. Since $O_{k_1}\bigcup\lt(\bigcup_{j=2}^{m}T_j\rt)$ is connected and $\Omega_{k_0}\subset O_{k_1}$, by definition of $\Omega_k$, we have $O_{k_1}\bigcup\lt(\bigcup_{j=2}^{m}T_j\rt)\subset \Omega_k$. Since $\Omega=\bigcup O_k$, it follows that (3) is satisfied.
		
		Now we claim that each $\Omega_k$ is simply-connected. We argue by contradiction. Suppose $\Omega_k$ is not simply-connected for some $k$. Then \xred{we can find some} closed curve $\Gamma\subset \Omega_k$ such that there exists $x\in \mathrm{Int}(\Gamma)\cap (O_k)^c$. By the definition of $O_k$, we have $\mathrm{dist}(x,\partial\Omega)\leq 2^{-k}$. Let $y\in\partial\Omega$ be such that $|x-y|=\mathrm{dist}(x,\partial\Omega)$, and let $z$ be the intersection of $\Gamma$ with the line segment joining $x$ and $y$. Then clearly we have $|z-y|\xred{\leq}|x-y|\leq 2^{-k}$. On the other hand, since $z\in\Gamma\subset\Omega_k$, we have $|z-y|\geq \mathrm{dist}(z,\partial\Omega)>2^{-k}$. This is a contradiction. It follows that $\Omega_k$ is simply-connected.
		
		\em Step 2: proof of Lemma completed. \rm Without loss of generality we can assume $0\in \Omega_k$ for all $k$. For any $\ep\in (0,2^{-k})$, $v_{\ep}=v*\rho_{\ep}$ is such that $\mathrm{curl}v_{\ep}=0$ on $\Omega_k$. Since $\Omega_k$ is simply-connected, there exists $f_{\ep}$ such that $\na f_{\ep}=v_{\ep}$ on 
		$\Omega_k$ and $f_{\ep}(0)=0$. Now take some sequence $\ep_n\rightarrow 0$. By basic properties of convolutions, we know $\na f_{\ep_n}\overset{L^p(\Omega_k)}{\rightarrow} v$ as $n\rightarrow \infty$, for all $1\leq p<\infty$. 
		
		Since $v\in L^{\infty}(\Omega;\R^2)$, we have $\lVert v_{\ep}\rVert_{\infty}\leq \lVert v\rVert_{\infty}$, and hence $\{f_{\ep_n}\}$ is a sequence of equicontinuous functions on $\overline{\Omega}_k$ with $f_{\ep_n}(0)=0$. It follows from the Arzel\`{a}-Ascoli Theorem that for some subsequence (not relabeled) 
		$f_{\ep_n}\overset{L^{\infty}(\Omega_k)}{\rightarrow} f_k$ for some Lipschitz function $f_k$ with $f_k(0)=0$. Therefore $\na f_k=v$ a.e. on $\Omega_k$. 
		
		We claim 
		\begin{equation}
		\label{curleq1}
		f_l=f_k\text{ on }\Omega_k \text{ for all } l>k.
		\end{equation}
		Indeed, the equation (\ref{curleq1}) follows from the facts that $f_l-f_k$ is Lipschitz and $f_l(0)=f_k(0)$ and 
		$\na (f_l-f_k)=0$ a.e. on $\Omega_k$. Thus by (\ref{curleq1}) we can define 
		$$
		f(x)=\lt\{\begin{array}{ll}f_1(x)&\text{ on }\Omega_1\nn\\
		f_2(x)&\text{ on }\Omega_2\nn\\ 
		&\dots \nn\\
		f_k(x)&\text{ on }\Omega_k\nn\\
		&\dots \nn\\
		 \end{array}\nn\rt\}.
		$$
		And finally $\na f=v$ a.e. on $\Omega$.
		\end{proof}
				
		Finally, we provide the proof of Lemma \ref{l16}.
				
		\begin{proof}[Proof of Lemma \ref{l16}]
		We mostly follow the proof of Lemma 4 in \cite{mul2}. Let us consider the function $\varphi$ defined by 
		\begin{equation*} %\label{l16.1}
		\varphi(z)=
		\begin{cases}
		z\cdot \xi &\text{for } z\cdot \xi >0,\\
		0 &\text{for } z\cdot \xi \leq 0,
		\end{cases}
		\end{equation*}
		and the map $F$ given by 
		\begin{equation*}
		F(z)=
		\begin{cases}
		\xi &\text{for } z\cdot \xi >0,\\
		0 &\text{for } z\cdot \xi \leq 0.
		\end{cases}
		\end{equation*}
		Note that $F$ is the gradient of $\varphi$ whenever $\varphi$ is differentiable.
					
		Now we construct a sequence $\{\varphi_k\}_{k}$ in $C_{c}^{\infty}(\Rb^2)$ such that
		\begin{equation}\label{l16.2}
		\{(\varphi_k(z),\nabla\varphi_k(z))\}_{k} \qd \text{ is bounded uniformly for bounded } z,
		\end{equation}
		\begin{equation}\label{l16.3}
		(\varphi_k(z),\nabla\varphi_k(z)) \overset{k\rightarrow \infty}{\rightarrow} (\varphi(z),F(z)) \qd\text{ for all } z.
		\end{equation}
		Here we use an approximation that was used by the first author in \cite{Lo} to make the proof more transparent than that in \cite{mul2}. Clearly there exists a monotone smooth function $s_0: \Rb\rightarrow\Rb$ such that $s_0(x)\equiv 0$ for $x\leq 0$ and $s_0(x)=x$ for $x\geq 1$. Given $k\in \mathbb{N}^{+}$, define $s_k(x):= \frac{1}{k}s_0(kx)$. It is easy to check that $s_k$ is a smooth function satisfying 
		\begin{equation}\label{l16.6}
		\{(s_k(x), s_k'(x)) \} \qd \text{ is bounded uniformly for bounded } x,
		\end{equation}
		\begin{equation}\label{l16.7}
		(s_k(x), s_k'(x)) \overset{k\rightarrow \infty}{\rightarrow} (s(x), f(x)) \qd\text{ for all } x, 
		\end{equation}
		where
		\begin{equation*}
		s(x)=
		\begin{cases}
		x & \text{ for } x>0,\\
		0 & \text{ for } x\leq 0,
		\end{cases}
		\end{equation*}
		and 
		\begin{equation*}
		f(x)=
		\begin{cases}
		1 & \text{ for } x>0,\\
		0 & \text{ for } x\leq 0.
		\end{cases}
		\end{equation*}
		Now we define $\varphi_k(z)=s_k(z\cdot \xi)\chi_k$, where $\chi_k\in C_{c}^{\infty}(\Rb^2)$ satisfies $\spt(\chi_k)\subset\subset B_{k+1}(0)$ and $\chi_k\equiv 1$ on $B_k(0)$. It is clear that $\varphi_k\in C_{c}^{\infty}(\Rb^2)$ and $\nabla \varphi_k(z)=s_k'(z\cdot \xi)\xi$ for $z\in B_k(0)$. One can check directly that the properties \eqref{l16.6}-\eqref{l16.7} for $s_k$ translate to \eqref{l16.2}-\eqref{l16.3}. 
					
		According to Lemma \ref{l13},
		\begin{equation*}
			\Phi_k(z):= \varphi_k(z)z + \lt(\nabla\varphi_k(z)\cdot z^{\perp}\rt)z^{\perp}
		\end{equation*}
		is an entropy. It is clear that \eqref{l16.2} implies that $\{\Phi_{k}(z)\}$ is bounded uniformly for bounded $z$. According to \eqref{l16.3}, 
		\begin{equation*}
			\Phi_k(z)\rightarrow \varphi(z)z + \lt(F(z)\cdot z^{\perp}\rt)z^{\perp}=
			\begin{cases}
			|z|^2 \xi &\text{for } z\cdot \xi>0,\\
			0 &\text{for } z\cdot \xi \leq 0,
			\end{cases}
		\end{equation*}
		which is \eqref{l16.5}.
		\end{proof}

\end{document}